\documentclass[11pt]{article}

\usepackage{geometry}
\usepackage[pdftex]{graphicx}
\usepackage{amsfonts,amssymb,amsmath,amsthm}
\usepackage{epstopdf}
\usepackage{longtable}
\usepackage[normalem]{ulem}
\usepackage{graphicx,multirow,float}
\usepackage{hyperref,verbatim,float,stackrel}

\newcommand*{\cB}{\mathcal{B}}
\newcommand*{\cP}{\mathcal{P}}
\newcommand*{\cS}{\mathcal{S}}
\newcommand*{\bbC}{\mathbb{C}}
\newcommand*{\bbI}{\mathbb{I}}
\newcommand*{\bbN}{\mathbb{N}}
\newcommand*{\bbR}{\mathbb{R}}
\newcommand*{\bbP}{\mathbb{P}}
\newcommand*{\bbZ}{\mathbb{Z}}
\newcommand*{\rB}{\mathrm{B}}
\newcommand*{\rH}{\mathrm{H}}
\newcommand*{\rL}{\mathrm{L}}

\newcommand{\be}[1]{\begin{equation} #1 \end{equation}}
\newcommand{\bes}[1]{\begin{equation*} #1 \end{equation*}}
\newcommand{\ea}[1]{\begin{align} #1 \end{align}}
\newcommand{\eas}[1]{\begin{align*} #1 \end{align*}}
\numberwithin{equation}{section}

\newtheorem{theorem}{Theorem}
\newtheorem{remark}[theorem]{Remark}
\newtheorem{example}[theorem]{Example}
\newtheorem{lemma}[theorem]{Lemma}
\newtheorem{corollary}[theorem]{Corollary}
\newtheorem{proposition}[theorem]{Proposition}
\newtheorem{definition}[theorem]{Definition}
\newtheorem{problem}[theorem]{Problem}

\numberwithin{theorem}{section}

\newcommand{\thm}[1]{\begin{theorem} #1 \end{theorem}}

\newcommand{\prop}[1]{\begin{proposition} #1 \end{proposition}}
\newcommand{\rem}[1]{\begin{remark} #1 \end{remark}}
\newcommand{\lem}[1]{\begin{lemma} #1 \end{lemma}}
\newcommand{\cor}[1]{\begin{corollary} #1 \end{corollary}}

\newcommand{\prf}[1]{\begin{proof} #1 \end{proof}}

\def\R#1{(\ref{#1})}
\newcommand*{\nm}[1]{{\left\|#1\right \|}} 
\newcommand{\abs}[1]{\left|#1\right|}
\newcommand*{\ip}[2]{\langle #1, #2 \rangle}
\newcommand{\keywords}[1]{\textbf{{Keywords:}} \small{#1}}
\newcommand{\AMS}[1]{\textbf{{AMS:}} \small{#1}}

\def\D{\,\mathrm{d}}
\def\I{\mathrm{i}}
\def\E{\mathrm{e}}
\newcommand*{\rth}{\mathrm{th}}
\newcommand*{\supp}{\mathrm{supp}}

\title{Density theorems for nonuniform sampling of bandlimited functions using derivatives or bunched  measurements} 

\author{Ben Adcock \thanks{Department of Mathematics, Simon Fraser University, Canada (ben\_adcock@sfu.ca).} \and Milana Gataric \thanks{DPMMS, Centre for Mathematical Sciences, University of Cambridge, UK (m.gataric@maths.cam.ac.uk).} \thanks{The corresponding author.}  \and Anders C. Hansen \thanks{DAMTP, Centre for Mathematical Sciences, University of Cambridge, UK; Department of Mathematics, University of Oslo, Norway (ach70@cam.ac.uk).} 
}

\begin{document}
\maketitle

\begin{abstract}
We provide sufficient density condition for a set of nonuniform samples to give rise to a set of sampling for multivariate bandlimited functions when the measurements consist of pointwise evaluations of a function and its first $k$ derivatives. Along with explicit estimates of corresponding frame bounds, we derive the explicit density bound and show that, as $k$ increases, it grows linearly in $k+1$ with the constant of proportionality $1/\E$. Seeking larger gap conditions, we also prove a multivariate perturbation result for nonuniform samples that are sufficiently close to sets of sampling, e.g.\  to uniform samples taken at $k+1$ times the Nyquist rate.

Additionally, in the univariate setting, we consider a related problem of so-called nonuniform bunched sampling, where in each sampling interval $s+1$ bunched measurements  of a function are taken and the sampling intervals are permitted to be of different length. We derive an explicit density condition which grows linearly in $s+1$ for large $s$, with the constant of proportionality depending on the width of the bunches. The width of the bunches is allowed to be arbitrarily small, and moreover, for sufficiently narrow bunches and sufficiently large $s$, we obtain the same result as in the case of univariate sampling with $s$ derivatives.
\end{abstract}


\vspace{3mm}

\AMS{42C15, 94A20, 41A05}

\section{Introduction}

In this paper we consider two related sampling scenarios. The first one assumes that a bandlimited function $f$ and its first $k$ derivatives are sampled at nonuniformly spaced points. This problem is motivated by applications in seismology, where certain recently-developed detectors are able to measure both $f$ and its spatial gradient. However, there are also various other applications and for different examples we direct the reader to \cite{EldarFilterbank, Lazaro} and  references therein. The second scenario considered in this paper assumes that, instead of evaluating derivatives of $f$ at $\{ x_n \}_{n \in I}$, $f$ is measured at an additional $s$ sampling points around each $x_n$. One can think of this scenario as function $f$ being evaluated at $s+1$ different nonuniform sampling sequences. When these sequences are uniform, the problem is known as bunched sampling or recurrent nonuniform sampling and has been extensively studied in literature (for more details see below). 

The purpose of this paper is to understand the gain one can expect by nonuniformly sampling derivatives or by nonuniformly sampling at bunched points. We derive explicit sufficient conditions for stable recovery in terms of densities of sampling points. In particular, we show that the maximal allowed gap between sampling points (or bunches of sampling points) grows linearly in $k+1$ (or $s+1$) for large $k$ (or $s$), which translates into increasing savings in cost and effort.

\subsection{Nonuniform sampling of bandlimited functions}
The topic of nonuniform sampling has been extensively researched in the last several decades. See, for example, \cite{AldroubiGrochenigSIREV,BenedettoBook, BenedettoSpiral,FeichtingerGrochenigIrregular, GrochenigStrohmerMarvasti,Marvasti,SeipBook,young} and references therein.  Let $\{ x_n \}_{n \in I}$ be a set of points in $\bbR^d$.  Then $\{ x_n \}_{n \in I}$ is a \textit{set of sampling} if there exist  $A,B>0$ such that
\be{
\label{set_of_samp}
A \| f \|^2 \leq \sum_{n \in I} | f(x_n) |^2 \leq B \| f \|^2,
}
for all bandlimited functions $f \in \rB(\Omega) = \{ g \in \rL^2(\bbR^d) : \supp(\hat{g}) \subseteq \Omega \}$, where $\Omega$ is a compact set.  This condition is equivalent to the existence of particular frame for $\rB(\Omega)$.  When \R{set_of_samp} holds, stable reconstruction of $f$ from the samples $\{ f(x_n) \}_{n \in I}$ is possible, and this can be carried out via a number of different algorithms.  Note that in practice, \R{set_of_samp} is usually replaced by a weighted sum, to accommodate clustering of the sampling points \cite{ AldroubiGrochenigSIREV,GrochenigIrregular,GrochenigIrregularExpType,GrochenigStrohmerMarvasti} and recently \cite{BAMGACHNonuniform1D,AGH2DNUGS}.

Since \R{set_of_samp} implies stable recovery, it is important to provide sufficient conditions for \R{set_of_samp} to hold in terms of  \textit{density} of the set $\{ x_n \}_{n \in I}$.  In one dimension, an (almost) full characterization is known in terms of the so-called Beurling density, due to  Landau \cite{Landau}, Jaffard \cite{Jaffard} and Seip \cite{SeipJFA} (see also \cite{ChristensenFramesAMS}).  In higher dimensions a result of Beurling \cite{BeurlingDiffOp} provides a sharp, sufficient condition whenever the points $\{ x_n \}_{n \in I}$ are relatively separated and $\Omega$ is the unit Euclidean ball.  This was extended to any compact, convex and symmetric $\Omega$ in \cite{BenedettoSpiral} and \cite{OlevskiiUlanovskii}.  Beurling's result has also been shown to hold when weights are incorporated into \R{set_of_samp} to compensate for clustering \cite{AGH2DNUGS}. Unfortunately, such results do not give explicit estimates for the frame bounds $A$ and $B$, and in particular, their ratio $\sqrt{B/A}$, which determines the stability of any reconstruction algorithm.  Explicit estimates  were provided by Gr\"ochenig \cite{GrochenigIrregular}.  In one dimension, Gr\"ochenig's result is sharp, but in higher dimensions, the density condition deteriorates with $d$ and thus ceases to be sharp.  In recent work \cite{AGH2DNUGS}, the dimension dependence of such conditions has been improved.  For example, when $\Omega$ is the unit Euclidean ball, a result of \cite{AGH2DNUGS} gives explicit frame bounds for a density condition that---although somewhat more restrictive than Beurling's---is independent of $d$.

Another approach to obtain guarantees for \R{set_of_samp} which allows for larger gaps between sampling points is to prove perturbation theorems.  That is, if a set $\{ x_n \}_{n \in I}$ is known to give a set of sampling, and if $\{ \tilde{x}_n \}_{n \in I}$, then one seeks to establish conditions on the maximal distance $\epsilon = \sup_{n \in I} | \tilde{x}_n - x_n |$  such that $\{ \tilde{x}_n \}_{n \in I}$ also satisfies \R{set_of_samp} with explicit bounds $\tilde{A}$ and $\tilde{B}$ depending on $A$, $B$ and $\epsilon$.  In particular, since uniform samples taken at the Nyquist rate provide sets of sampling with $B/A = 1$, perturbation theorems allow for nonuniform sampling points to be taken with larger separation, provided they are sufficiently close to uniform sampling points.  In one dimension, Kadec's-$1/4$ theorem \cite{young} (see also \cite{BalanFourierFramesKadec,ChristensenFramesAMS}) is a sharp perturbation theorem in the sense that condition on $\epsilon$ cannot be increased.  Higher-dimensional generalizations have been obtained in \cite{BaileySampRecov,ChuiShuiTrigPerturb,FavierZalekFramesPerturb,SunZhouMultivarTrigPerturb}.

\subsection{Derivatives sampling}

Uniform sampling of derivatives is a classical topic in sampling theory, see \cite{FogelSampling,JagermanFogel,LindenAbramson,PapoulisGenSamp,RawnDerivatives,ZibulskiPeriodicDerivatives} and references therein.  It is known that one can exceed the Nyquist criterion by a factor of $k+1$ by sampling $f$ and its first $k$ derivatives \cite{PapoulisSignalAnalysis,RawnDerivatives}.  On the other hand, relatively few papers have considered nonuniform sampling with derivatives.  In the univariate setting, by extending Gr\"ochenig's results \cite{GrochenigIrregular} for univariate nonuniform sampling to the case of derivatives, it was shown in \cite{RazafinjatovoDerivs} that the maximum allowable spacing between sampling points grows asymptotically as a linear function of $k+1$, with constant of proportionality equal to $1/\E$ (Gr\"ochenig also considered the case of derivatives in  \cite{GrochenigIrregular}, but did not show the improvement in the maximum spacing).  Also, in the univariate setting, periodic nonuniform sampling with derivatives was considered in \cite{ZibulskiPeriodicDerivatives}. Multivariate nonuniform  sampling with derivatives was addressed in \cite{LandauConditionsGroch}, where necessary density conditions are derived.  

To the best of our knowledge, no work has addressed sufficient guarantees for stable sampling with derivatives in the multivariate setting. This is the main task of the present paper.

\subsection{Bunched sampling}

Uniform bunched sampling---also known as recurrent nonuniform sampling since it assumes periodic groups of nonuniform samples---has been a topic of numerous  papers, in both the univariate \cite{BunchedSampButzer,EldarFilterbank,BunchedSampKohlenberg,PapoulisGenSamp,PapoulisSignalAnalysis,SommenRecurrent,BunchedSampStrohmer,BunchedSampYen} and the multivariate case \cite{FaridaniBunchedSamp,MultidimFilterbank}. In this setting, as in the uniform derivative sampling, one is allowed to sample above the Nyquist rate. Namely, if the uniform  set of bunched points is interpreted as the union of $s+1$ uniform sequences, then each of these sequences can be taken at $s+1$ times the Nyquist rate. However, one might want to know what happens if 
each of these sequences is nonuniform, i.e.~if the groups of
nonuniform samples are not repeating periodically, but instead are distributed nonuniformly. This setting corresponds to nonuniform bunched sampling that we consider in the present paper. 

To the best of our knowledge, there is no earlier work which considers nonuniform bunched sampling. Nonetheless, there is a list of questions that naturally arises in regards to this sampling problem. For instance, in the setting where one must allow for bigger distances between sampling sensors due to some natural constraints, could one use nonperiodic bunches of sensors and thus sample above the Nyquist rate? Is there a lower bound on the width of those bunches? We answer these questions in the univariate setting in the current paper.

\subsection{Function recovery}

This paper focuses on the question of when a stable set of sampling exists, rather
than on how to recover $f$ from its samples. For practical applications, it is also important that one  considers the reconstruction problem, where only finitely-many samples are given. Using techniques of generalized sampling \cite{BAMGACHNonuniform1D,AGH2DNUGS,BAACHShannon,BAACHOptimality}, it is possible to show that stable reconstruction is possible in any given finite-dimensional approximation subspace, subject to two conditions on the data. First, it should satisfy the same density requirement that ensures a stable set of sampling. Second, it should cover a sufficiently large region of $\bbR^d$, with the precise size of this region depending solely on the particular choice of the approximation subspace. 

In particular, this applies to the problem of recovering a bandlimited function from its pointwise samples and samples of its derivatives, or its pointwise samples at bunched points. In this case, one retains the same density requirements as derived in the current paper. Additionally, one would need to choose a suitable approximation space for recovery of a bandlimited function. For this purpose, one could choose the finite-dimensional space of trigonometric polynomials, for example, as it is the case in the well-known ACT algorithm \cite{FeichtingerEtAlEfficientNonuniform,GrochenigIrregularExpType,GrochenigModernSamplingBook,GrochenigStrohmerMarvasti}. By this choice of the reconstruction space, one extends the ACT algorithm in a straightforward manner to the case of derivative or bunched sampling, with the density guarantee that increases linearly in $k+1$ or $s+1$, respectively. See also \cite{RazafinjatovoRec}. 

For additional literature on the problem of recovering an element of a Hilbert space from its samples, see for example \cite{GenSampProbability,eldar2003FAA,eldar2005general,hrycakIPRM,unser1994general}.

\subsection{Main results of the paper}

The major part of the paper, Section \ref{s:derivatives}, is dedicated to derivatives sampling. In particular, \S \ref{ss:multivariate_thms} deals with the multivariate derivatives sampling, while \S \ref{ss:bound_1D} and \S \ref{ss:bound_tensor} address special cases of univariate and line-by-line sampling. Additionally, in \S \ref{s:perturbation}, a perturbation result for derivatives sampling is given. In Section \ref{s:bunched_samp},  one-dimensional nonuniform bunched sampling is considered, in context of fusion frames \S \ref{ss:bunched_fusion} and regular frames \S \ref{ss:bunched_frames}.

Our first result, Theorem \ref{t:main_frame_dD}, provides an upper bound on the maximum allowable sampling density $\delta$ (see \R{delta_def} for a definition), such that samples of derivatives give rise to a particular frame. The density bound as well as the explicit estimates of the corresponding frame bounds depend on the number of derivatives $k$, the norm used in specifying $\delta$ and a certain geometric property of the domain $\Omega$.  For large $k$, the maximum allowed $\delta$ grows linearly in $k+1$ with constant of proportionality $1/\E$.  This extends the univariate result of \cite{RazafinjatovoDerivs} to the multivariate setting, as well as the multivariate $k=0$ (no derivative) results of \cite{GrochenigIrregular,GrochenigModernSamplingBook}, and recently \cite{AGH2DNUGS}, to the case of derivatives.

In our second result, Theorem \ref{t:main_fram_1D}, we present a univariate density condition that leads to a small improvement over \cite{RazafinjatovoDerivs} for $k \geq 2$ derivatives.  This follows the technique of \cite{GrochenigIrregular} for the univariate case based on Wirtinger inequalities.  We provide an explicit calculation of the optimal constants in certain higher-order Wirtinger inequalities, which, replicating the techniques of \cite{GrochenigIrregular} for the case of derivatives, lead to modestly  improved estimates for $\delta$ for finite $k$.  Such improved bounds can be used to get better estimates for two-dimensional spatial-temporal sampling scenarios, as we consider in Proposition \ref{p:main_tensor_1}.

Next, we provide Theorem \ref{t:perturb}, which gives a perturbation estimate for nonuniform sampling with derivatives.  We show that if $\{ x_n \}_{n \in I}$ is a stable set of sampling for derivatives, then so is $\{ \tilde{x}_n \}_{n \in I}$ whenever $\sup_{n \in I} | x_n - \tilde{x}_n |$ is sufficiently small.  In particular, small perturbations of the uniform sampling points taken at $k$ times the Nyquist give rise to stable sets of sampling. This extends existing results given in \cite{BaileySampRecov, SunZhouMultivarTrigPerturb} to the case of sampling with derivatives. Moreover, it improves those results since we provide a dimension independent bound for appropriate domains $\Omega$.

In Section  \ref{s:bunched_samp}, we address univariate nonuniform bunched sampling and, in Theorem \ref{t:bunched_samp}, we give density guarantees in order to obtain a particular fusion frame \cite{CasazzaKutyniokFusion,CasazzaKutyniokLi}. Similarly as in derivatives sampling, we show that the density bound increases linearly with $s+1$ (the number of samples in each bunch) with constant of proportionality depending on the width of the bunches. In Proposition \ref{p:bunched_samp_finitediff}, we derive the same density condition under which it is possible to construct a particular frame based on divided differences. The  points within the same bunch are permitted to get arbitrarily close to each other, since we use appropriate weights.  Furthermore, we show that the corresponding density bound in the limit---for small width of bunches and for large $s$---gives the same density bound as the one we provide for the univariate sampling with $s$ derivatives, i.e the density bound grows linearly in $s+1$ with the constant of proportionality $1/\E$.

\section{Preliminaries}\label{s:pre}
Let $d \geq 1$ denote dimension.  For $x = (x_1,\ldots,x_d) \in \bbR^d$ and $p \geq 1$, we write $| x |_p$ for the $\ell^p$-norm, i.e.\ $| x |_p = \left ( \sum^{d}_{j=1} | x_j |^p \right )^{1/p}$, when $1\leq p<\infty$, and $|x|_{\infty}=\max_{j=1,\ldots,d} |x_j|$.  If $y = (y_1,\ldots,y_d) \in \bbR^d$ we write $x \cdot y = x_1 y_1 + \ldots + x_d y_d$ for the dot product of $x$ and $y$.  We let $\rL^2(\bbR^d)$ be the space of square-integrable functions on $\bbR^d$ with inner product
\bes{
\ip{f}{g} = \int_{\bbR^d} f(x) \overline{g(x)} \D x,
}
and corresponding norm $\| f \| = \sqrt{\ip{f}{f}}$.  We denote the Fourier transform by
\bes{
\hat{f}(\omega) = \int_{\bbR^d} f(x) \E^{-\I \omega \cdot x } \D x,\quad \omega \in \bbR^d.
}
Note that Parseval's identity reads $\| \hat{f} \|^2 = (2 \pi)^d \| f \|^2$ with this normalization.

Throughout the paper we write $D^{\alpha}= \partial^{|\alpha|_1}/\partial^{\alpha_1}_{x_1} \cdots \partial^{\alpha_d}_{x_d}$ for the partial derivative operator.
Here $\alpha = (\alpha_1,\ldots,\alpha_d) \in \bbN^d_0$, $\bbN_0 = \bbN \cup \{0\}$, and $| \alpha |_1 = \alpha_1+\ldots+\alpha_d$.  If $\alpha ! = \alpha_1 ! \alpha_2 ! \cdots \alpha_{d} !$ then the multinomial formula reads
\be{
\label{multinomial}
\left ( \sum^{d}_{j=1} x_j \right )^k = \sum_{| \alpha |_1 = k} \frac{k!}{\alpha!} x^\alpha,\quad k \in \bbN, x \in \bbR^d,
}
where $x^{\alpha}$ is given by $x^{\alpha} = \prod^{d}_{j=1} (x_j)^{\alpha_j}$.

Let $\Omega \subseteq \bbR^d$ be a compact set.  The space of $\Omega$-bandlimited functions is defined by
\bes{
\rB(\Omega) = \{ f \in \rL^2(\bbR^d) : \mathrm{supp}(\hat{f}) \subseteq \Omega \}.
}
Observe also that $\rB(\Omega)$ is closed with respect to differentiation, and moreover we have the multivariate Bernstein inequality
\be{
\label{Bernstein}
\| D^{\alpha} f \| \leq \bar{\omega}^{\alpha} \| f \|,\quad \forall f \in \rB(\Omega),
}
where $\bar{\omega} = (\bar{\omega}_1,\ldots,\bar{\omega}_d)^{\top}$ and $\bar{\omega}_j = \sup_{\omega \in \Omega} | \omega_j |$.
In this paper, we measure the density of the sampling points in terms of an arbitrary norm $\abs{\cdot}_{*}$ on $\bbR^d$. Beside standard $\ell^p$-norms, one might want to use the norm induced by the polar set of $\Omega$, for example, cf.~Remark \ref{r:delta_lit}. Since all norms are equivalent on $\bbR^d$, we let $a,b>0$ be the optimal constants such that
\be{
\label{ab_star_def}
a | x |_{*}  \leq | x|_2 \leq b | x |_{*},\quad \forall x \in \bbR^d.
}
We shall also define the quantity
\be{
\label{m_Omega_def}
m_{\Omega} = \sup_{\omega \in \Omega} | \omega |_2.
}
Finally, let us now recap the notion of frames \cite{christensen2003introduction}.  Let $\rH$ be a Hilbert space with norm $\nm{\cdot}$ and inner product $\ip{\cdot}{\cdot}$.  A set $\{ \phi_n \}_{n \in I} \subseteq \rH$, where $I$ is an index set, is called a frame for $\rH$ if there exist $A,B>0$ such that
\bes{
A \| f \|^2 \leq \sum_{n \in I} | \ip{f}{\phi_n} |^2 \leq B \| f \|^2,\quad \forall f \in \rH.
}
We refer to $A$ and $B$ as the \textit{frame bounds}.

\section{Nonuniform sampling with derivatives}\label{s:derivatives}

Let $\{ x_n \}_{n \in I} \subseteq \bbR^d$ be a set of sampling points, where $I$ is a countable index set.  Let $f \in \rB(\Omega)$, and suppose that we are given the measurements
\bes{
D^{\alpha} f(x_n),\quad n \in I, | \alpha |_1 \leq k.
}
As discussed, stable recovery is possible if there exist nonnegative \textit{weights} $\mu_{n,\alpha}$ such that
\be{
\label{samples_sum_stab}
A \| f \|^2 \leq \sum_{n \in I} \sum_{|\alpha|_1 \leq k} \mu_{n,\alpha} | D^{\alpha} f(x_n) |^2 \leq B \| f \|^2,\quad \forall f \in \rB(\Omega),
}
holds for constants $A,B>0$. Following \cite{RazafinjatovoDerivs}, let us now define the function
\be{\label{Phi_fun}
\Phi_{\Omega}(x) =  \frac{1}{(2 \pi)^d} \int_{\Omega} \E^{\I \omega \cdot x} \D \omega,\quad x \in \bbR^d.
}
For a given $f \in \rB(\Omega)$, we have $\hat{f}(\omega) = \hat{f}(\omega) \bbI_{\Omega}(\omega)$. If $g\in \rB(\Omega)$ is such that $\hat{g}(\omega) = \bbI_{\Omega}(\omega)$, then by the convolution theorem we can write 
\bes{
f(x)  =  \int_{\bbR^d} f(s) g(x-s) \D s = \frac{1}{(2 \pi)^d} \int_{\bbR^d} \int_{\Omega} f(s) \E^{\I \omega \cdot (x-s) }  \D \omega \D s = \ip{f}{\Phi_{\Omega}(\cdot-x) }.
}
Therefore 
\bes{
D^{\alpha} f(x_n) = \ip{D^\alpha f}{\Phi_\Omega(\cdot-x_n) } = (-1)^{|\alpha|_1}\ip{f}{D^\alpha \Phi_\Omega(\cdot-x_n)}.
}
Hence \R{samples_sum_stab} is equivalent to the condition that the set of functions
\bes{
\left \{ \sqrt{\mu_{n,\alpha}}D^{\alpha} \Phi_\Omega(\cdot-x_n) : n \in I, |\alpha|_1 \leq k \right \},
}
forms a frame for $\rB(\Omega)$ with frame bounds $A,B>0$.

Similarly, after differentiation and using Parseval's identity, \R{samples_sum_stab} becomes 
\bes{
{A}/{(2\pi)^d} \ \| \hat f \|^2 \leq \sum_{n \in I} \sum_{|\alpha|_1 \leq k} \mu_{n,\alpha} | \ip{\hat f}{(-\I\omega)^{\alpha}\E^{-\I\omega x_n}} |^2 \leq {B}/{(2\pi)^d}\ \| \hat f \|^2,\quad \forall f \in \rB(\Omega)
}
and therefore,  \R{samples_sum_stab} is also equivalent to $\left\{(2\pi)^{d/2}\sqrt{\mu_{n,\alpha}}(-\I\omega)^{\alpha}\E^{-\I\omega x_n}: n\in I,|\alpha|_1 \leq k \right\}$
being a Fourier frame for $\rL^2(\Omega)= \{ f \in \rL^2(\bbR^d) : \mathrm{supp}(f) \subseteq \Omega \}$ with the frame bounds $A,B>0$; see, for example, \cite{young}.

In what follows, we provide sufficient conditions for \R{samples_sum_stab} for an appropriate choice of weights.  As is standard in nonuniform sampling, our weights shall be related to the Voronoi cells $\{ V_n \}_{n \in I}$ of the sampling points $\{ x_n \}_{n \in I}$  \cite{FeichtingerGrochenigIrregular,GrochenigIrregular,GrochenigModernSamplingBook,GrochenigStrohmerMarvasti}.  Let $\abs{\cdot}_*$ be a norm on $\bbR^d$.  Given $n \in I$, we let
\bes{
V_n = \left \{ x \in \bbR^d : | x - x_n |_* \leq | x-x_m |_*,\  \forall m \in I, m \neq n \right \},\quad n \in I,
}
and define 
\be{
\label{weights}
\mu_{n,\alpha} = \frac{1}{\alpha!} \int_{V_n} (x-x_n)^{2 \alpha} \D x,\quad \alpha \in \bbN^d_0,\  n \in I.
}
As is also standard, our sufficient conditions for \R{samples_sum_stab} with these weights will be in terms of the density of the sampling points, measured in the following sense:
\be{
\label{delta_def}
\delta = \sup_{x \in \bbR^d} \inf_{n \in I} | x-x_n |_*.
}
Our aim is to find the maximal allowable density $\delta$ for which \R{samples_sum_stab} holds.

\subsection{The multivariate case}\label{ss:multivariate_thms}
For $z \in (0,\infty)$,  define
\ea{
R_k(z) &= \exp(z) - \sum^{k}_{r=0} \frac{1}{r!} z^r, \label{R_k} \\
\sigma^*_d(z) &= \frac{z + \sqrt{z(d+z)}}{d}, \label{sigma_star} \\
h_k(z) &= \exp(z) R_k(z), \label{h_k}\\
g_{k,d}(z) &=  \left(1+2\sigma^*_d(z)\right)^{d/2} \exp(z/\sigma^*_d(z)) R_k(z). \label{g_k,d}
}
Note that both $h_{k}$ and $g_{k,d}$ have limiting value $0$ as $z \rightarrow 0^+$, and both increase monotonically to infinity as $z \rightarrow \infty$. Hence they have well-defined inverse functions $H_{k}(w)$ and $G_{k,d}(w)$ for $w \in (0,\infty)$.

Our main result is now as follows:

\thm{
\label{t:main_frame_dD}
Suppose that the weights $\mu_{n,\alpha}$ are given by \textnormal{\R{weights}} and let $\delta$ be as in \textnormal{\R{delta_def}}.  If
\be{
\label{delta_condition}
\delta < \frac{C(k,d)}{m_{\Omega} b},\qquad C(k,d) = \max\left\{H_k(1),G_{k,d}(1)\right\},
}
where $b$ and $m_\Omega$  are as in \textnormal{\R{ab_star_def}} and \textnormal{\R{m_Omega_def}} respectively, then 
\bes{
A \| f \|^2 \leq \sum_{n \in I} \sum_{| \alpha |_1 \leq k} \mu_{n,\alpha} | D^{\alpha}f(x_n) |^2 \leq B \| f \|^2,\quad \forall f \in \rB(\Omega),
}
where $A,B>0$ satisfy
\be{\label{frame_bounds_thm}
A \geq \E^{-d} \left ( 1 - \min\left\{h_k(m_{\Omega} b \delta),g_{k,d}(m_{\Omega} b \delta)\right\}  \right )^2,\quad B \leq  \exp(2 m_{\Omega} b \delta + (m_{\Omega} b \delta)^2).
}
Equivalently, the set $\{ \sqrt{\mu_{n,\alpha}} D^{\alpha} \Phi_{\Omega}(\cdot -x_n) : n \in I, | \alpha |_1 \leq k \}$ forms a frame for $\rB(\Omega)$ with bounds $A$ and $B$.
}

The key part of this theorem---whose proof we defer to \S \ref{sss:main_thm_proof}---is the condition \R{delta_condition}. Note that an interesting facet of \R{delta_condition} is that it splits geometric terms depending on the domain $\Omega$ (the constant $m_{\Omega}$) and the norm used (encapsulated by the term $b$), from the nongeometric constant $C(k,d)$.  

Whilst values of $C(k,d)$ for fixed $k$ and $d$ are easily calculated and are presented in Table \ref{tab:delta_bound_dD}, to understand its behaviour it is interesting to consider the following two asymptotic regimes:
\bes{
\mbox{(i) $k$ fixed, $d \rightarrow \infty$,}\qquad \qquad  \mbox{(ii) $d$ fixed, $k \rightarrow \infty$}.
}
In (i) it is desirable for \R{delta_condition} to not decrease with $d$, i.e.\ the density bound does not worsen with increasing dimension.  For (ii), we desire linear increase in the bound with $k$, i.e.\ adding derivatives samples means that sampling points can be taken further apart, at as fast a rate as possible.

\begin{table}
\scalebox{0.82}{
\begin{tabular}{|c|c|c|c|c|c|c|c|c|c|c|c|c|c|c|c|c|c|c}
\hline
$k$ & 0 & 1 & 2 & 3 & 4 & $\ldots$ & 8 & 9 & $\ldots$ & 13 & 14 & $\ldots$ \\ \hline \hline
 $C(k,1)$  & 0.4812 & 0.8141 & 1.1268 & 1.4304 & \textit{1.7890} & $\ldots$ & \textit{3.2501} & \textit{3.6163} & $\ldots$ & \textit{5.0828}  & \textit{5.4498} & $\ldots$  \\ \hline 
 $C(k,2)$  & 0.4812 & 0.8141 & 1.1268 & 1.4304 & 1.7290 & $\ldots$ & 2.8976 & \textit{3.2424} & $\ldots$ & \textit{4.6462}  &  \textit{5.0000} & $\ldots$  \\ \hline
 $C(k,3)$  &  0.4812 & 0.8141 & 1.1268 & 1.4304 & 1.7290 & $\ldots$ & 2.8976 & 3.1862 & $\ldots$ & 4.3327  & \textit{4.6679} & $\ldots$  \\ \hline 
 $C(k,4)$  &  0.4812 & 0.8141 & 1.1268 & 1.4304 & 1.7290 & $\ldots$ & 2.8976 & 3.1862 & $\ldots$ & 4.3327  & 4.6180 & $\ldots$  \\ \hline 
 $C(k,5)$  &  0.4812 & 0.8141 & 1.1268 & 1.4304 & 1.7290 & $\ldots$ & 2.8976 & 3.1862 & $\ldots$ & 4.3327  & 4.6180 & $\ldots$ \\ \hline 
\end{tabular}} 
\scalebox{0.812}{
\begin{tabular}{|c|c|c|c|c|c|c|c|c|c|c|c|c|}
\hline
$k$ & $\ldots$ & 17 & 18 & 19 & 20 & 21 & 22 & 23 & 24 & 25 & 26   \\ \hline \hline
 $C(k,1)$  & $\ldots$ & \textit{6.5512} & \textit{6.9184} & \textit{7.2857} & \textit{7.6531} &  \textit{8.0205} & \textit{8.3879} & \textit{8.7553} &  \textit{9.1228} & \textit{9.4903}  & \textit{9.8578}  \\ \hline 
 $C(k,2)$  & $\ldots$ & \textit{6.0660} & \textit{6.4227} & \textit{6.7799} & \textit{7.1376} & \textit{7.4958} & \textit{7.8544} & \textit{8.2134} & \textit{8.5728} & \textit{8.9325}  &  \textit{9.2925}  \\ \hline
$C(k,3)$  &  $\ldots$ & \textit{5.7002} & \textit{6.0466} & \textit{6.3940} & \textit{6.7424} & \textit{7.0916} & \textit{7.4415} & \textit{7.7922} & \textit{8.1435} & \textit{8.4955}  & \textit{8.8480}  \\ \hline 
 $C(k,4)$  &  $\ldots$ & 5.4715 & 5.7553 & \textit{6.0813} & \textit{6.4207} & \textit{6.7614} & \textit{7.1031} & \textit{7.4457} & \textit{7.7893} & \textit{8.1338} & \textit{8.4791}  \\ \hline 
 $C(k,5)$  &  $\ldots$ & 5.4715 & 5.7553 & 6.0389 & 6.3223 & 6.0654 & 6.8883 & 7.1711 & \textit{7.4879} & \textit{7.8252} & \textit{8.1636}  \\ \hline 
\end{tabular}}
\caption{The constant $C(k,d)$ in the multidimensional density bound \textnormal{\R{delta_condition}}.  Italics indicate when $C(k,d)=G_{k,d}(1)$, and otherwise $C(k,d)=H_{k}(1)$.}
\label{tab:delta_bound_dD}
\end{table}

\subsubsection{Case (i)}
As seen in Table \ref{tab:delta_bound_dD}, the constant $C(k,d)$ is independent of $d$ for large $d$ and fixed $k$.  This follows from \R{delta_condition}, where it is clear that for large $d$ the maximum is achieved by $H_k(1)$, which is dimension independent, as opposed to $G_{k,d}(1)$ (it can be easily proved that $G_{k,d}(1)$ decreases with $d$).  Hence for large $d$, the only possible dimension-dependence in \R{delta_condition} arises from the factors $m_{\Omega}$ and $b$, which are determined by the domain $\Omega$ and the norm $\abs{\cdot}_*$ respectively.  For simplicity, suppose that $\Omega = \cB_{p}$ is the unit $\ell^p$-ball, $p > 0$, and let $\abs{\cdot}_{*} = \abs{\cdot}_{q}$, $q \geq 1$, be the $\ell^q$-norm.  Then \R{delta_condition} reads
\bes{
\delta < \frac{C(k,d)}{\max\{1,d^{1/2-1/p}\} \max\{1,d^{1/2-1/q}\}}.
}
This follows from the well-known inequality
\be{\label{norms_ineq}
\forall x \in \bbR^d, \quad |x|_{q}\leq|x|_{r}\leq d^{1/r-1/q} |x|_{q},\quad 0 < r \leq q.
}
In particular, if $p = q =2$ for example (i.e.\ $\Omega$ is contained in the unit Euclidean ball and the $\delta$-density is measured in the Euclidean metric), then \R{delta_condition} reduces to  $\delta < C(k,d)$.  For sufficiently large $d$, one therefore obtains the dimensionless bound $\delta < H_k(1)$.  On the other hand, if $\Omega = [-1,1]^d$ is the unit cube and $\abs{\cdot}_{*} = \abs{\cdot}_{2}$ is the $\ell^2$-norm, then we get square-root decay of the corresponding bound, which reads $\delta < H_k(1) / \sqrt{d}$ for large $d$.

\rem{
\label{r:k0_case}
The splitting of the bound \textnormal{\R{delta_condition}} into the factors $C(k,d)$ and $m_{\Omega} b$ is an extension of the result found in \textnormal{\cite{AGH2DNUGS}} to the case $k \geq 1$.  Therein the case $k=0$ was considered and the bound $\delta < \ln(2) / (m_{\Omega} b)$ was established.  Conversely, \R{delta_condition} reduces to the somewhat stricter condition $\delta < \ln \left ( \frac{1}{2}(1+\sqrt{5}) \right ) / (m_{\Omega} b)$ when $k=0$; note that $ \ln \left ( \frac{1}{2}(1+\sqrt{5}) \right ) \approx 0.4812 < 0.6931 \approx \ln(2)$.  This is due to the additional complications arising from a bound that holds for arbitrary many derivatives.}

\subsubsection{Case (ii)}
We now discuss the case of fixed $d$ and increasing $k$.  Empirically, Table \ref{tab:delta_bound_dD} and the left panel of Figure \ref{f:delta_bound_dD} show that, whilst $H_k(1)$ gives the better bound for small values of $k$, asymptotically for $k\rightarrow\infty$ the better bound is provided by $G_{k,d}(1)$.  We confirm this with the following lemma:

\begin{figure}
\begin{center}
\includegraphics[scale=0.61]{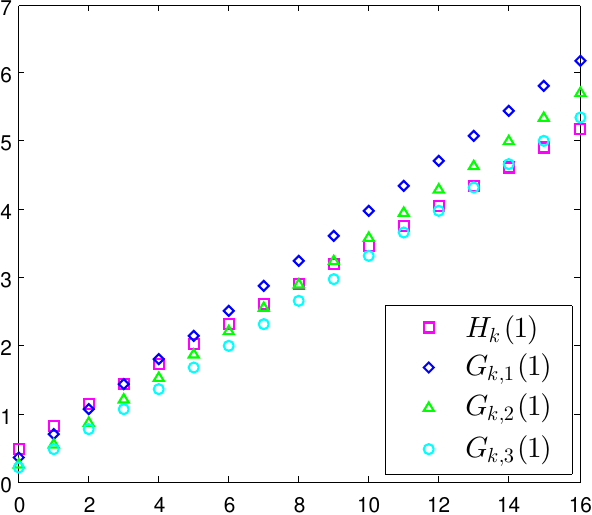} \quad
\includegraphics[scale=0.61]{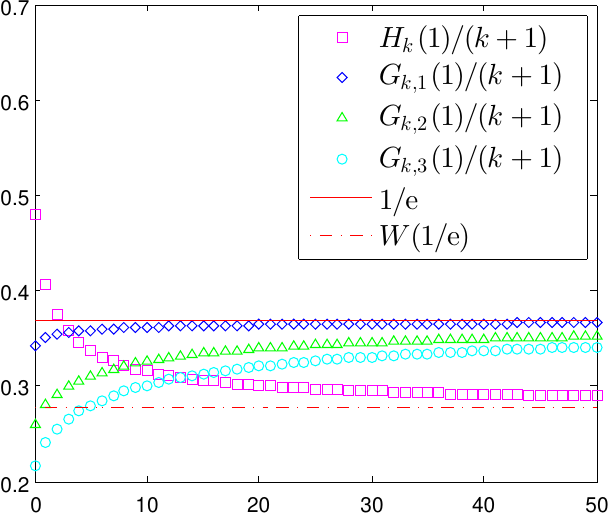} 
\caption{The constants in the multidimensional density bound \textnormal{\R{delta_condition}} $($left$)$ and their asymptotic behaviour $($right$)$.
} \label{f:delta_bound_dD}
\end{center}
\end{figure}

\lem{\label{l:large_k}
Let $W$ be the Lambert-$W$ function \textnormal{\cite{CorlessEtAlLambertW}}. We have
\begin{itemize}
 \item[$(a)$] $R_k(z)^{1/(k+1)} \sim \E z/(k+1)$ as $k \rightarrow \infty$, provided $z \leq c(k+1)$ for all large $k$ and some $c\in(0,1)$;
 \item[$(b)$] $H_k(1) \sim W\left({1}/{\E}\right) (k+1)$ as $k \rightarrow \infty$;
 \item[$(c)$] $G_{k,d}(1) \sim {1}/{\E} \ (k+1)$ as $k \rightarrow \infty$;
 \item[$(d)$] $C(k,d) \sim {1}/{\E} \ (k+1)$ as $k \rightarrow \infty$.
\end{itemize}
$($Note that $W\left({1}/{\E}\right)\approx0.2785$, while $1/\E\approx0.3679$.$)$
}

\prf{
To prove $(a)$, observe that 
\bes{
1 - \exp(-z) \sum^{k}_{r=0} \frac{1}{r!} z^r = \frac{\gamma(k+1,z)}{ \Gamma(k+1)} = P(k+1,z),
}
where $\gamma(\cdot,\cdot)$ is the lower incomplete Gamma function, and $\Gamma(\cdot)$ is the Gamma function \cite{AS}.  We  require an asymptotic expansion of $P(k+1,z)$ as $k \rightarrow \infty$ that is uniform in $z \leq c (k+1)$.  Such an expansion was obtained by Temme \cite{TemmeUniformGammaBeta,TemmeAsyGamma}.  It was shown there that
\bes{
P(\kappa,x) = \frac12 \mathrm{erfc}\left[-\eta (\kappa/2)^{1/2} \right ] - S_{\kappa}(\eta),
}
with
\bes{
S_{\kappa}(\eta) \sim (2 \pi \kappa)^{-1/2} \E^{-\kappa \eta^2/2} \sum^{\infty}_{k=0} c_k(\eta) \kappa^{-k},
}
as $\kappa \rightarrow \infty$, uniformly with respect to $\eta \in \bbR$, where $\kappa>0$, $x\geq0$, and
\bes{
\eta = \sqrt{2(\lambda-1 - \log(\lambda) ) },\quad \lambda = x/\kappa, \quad \mu = \lambda-1, 
}
with the square root having the same sign as $\mu$. Since $x/\kappa<1$, we have that $\mu<0$.  Here $\mathrm{erfc}(x) = 2 \pi^{-1/2} \int^{\infty}_{x} \E^{-t^2} \D t$ is the complementary error function and $c_k(\eta)$ are functions of $\eta$ only, with $c_0 = 1/\mu - 1/\eta$.
We require only the first term in the asymptotic expansion of $P(\kappa,x)$.   Since $\mathrm{erfc}(x) \sim \frac{\E^{-x^2}}{\sqrt{\pi} x}$ as $x \rightarrow \infty$, 
\bes{
P(\kappa,x) \sim -\frac{\exp(-\kappa \eta^2/2)}{\sqrt{2 \pi \kappa} \mu}=\frac{\exp(-\lambda \kappa + \kappa) \lambda^\kappa}{\sqrt{2 \pi \kappa} (1-\lambda) }, \quad \kappa\rightarrow \infty,
}
provided that $\lambda \leq c$ for some $c < 1$.  Set $k+1 = \kappa$ and $z= x =\kappa \lambda$.  Then, since $z \leq c(k+1)$ for some $c<1$, we get
\bes{
R_k(z)^{\frac{1}{k+1}} \sim \E \lambda = \E \frac{z}{k+1},
}
as $k \rightarrow \infty$ and $(a)$ follows.

To prove $(b)$, we shall use $(a)$.  Let $z = H_k(1)>0$, i.e.\ $h_k(z) = 1$.  We first show that there exists a $0 < c<1$ such that $z\leq c(k+1)$ for all large $k$.  Note that $R_k(z)\geq z^{k+1}/(k+1)!$.  Thus $z$ satisfies
\bes{
\exp(z) \frac{z^{k+1}}{(k+1)!} \leq 1.
}
Therefore
\bes{
\exp \left (\frac{z}{k+1} \right ) \frac{z}{k+1} \leq \frac{((k+1)!)^{\frac{1}{k+1}}}{k+1}.
}
By Stirling's formula, the right-hand side is asymptotic to $1/\E$ as $k \rightarrow \infty$.  Hence for large $k$, $z/(k+1) \leq W(1/\E) < 1$, as required.  We may now use $(a)$.  Since $h_k(z) = 1$ is equivalent to $\exp\left({\frac{z}{k+1}}\right) R_k(z)^{\frac{1}{k+1}}=1$, this now gives
\bes{
\exp\left({\frac{z}{k+1}}\right) \frac{z}{k+1} \sim \frac{1}{\E}, \quad k \rightarrow \infty.
}
Since the last identity is equivalent to
\bes{
\frac{z}{k+1} \sim W\left(\frac{1}{\E}\right), \quad k \rightarrow \infty,
}
we get the result.

We use a similar approach to prove $(c)$.  Let $z = G_{k,d}(1)$, i.e.\ $g_{k,d}(z) = 1$.  Then
\bes{
R_k(z) \leq \left(1+2\sigma^*_d(z)\right)^{d/2} \E^{z/\sigma^*_d(z)} R_k(z) = g_{k,d}(z) = 1,
}
and we deduce that $z \leq \frac{1}{\E} (k+1)$ as $k \rightarrow \infty$.  Hence we may apply $(a)$.  Note also that $z \rightarrow \infty$ as $k \rightarrow \infty$.  This follows from the fact that $\lim_{k\rightarrow\infty}g_{k,d}(z)=0$ for fixed $z$ and $d$.  Therefore, the equation $g_{k,d}(z) = 1$ can be written as
\bes{
\left(1+\frac{4z}{d}\right)^{\frac{d}{2(k+1)}} \exp\left(\frac{d}{2(k+1)}\right)  \frac{\E z}{k+1} \sim 1, \quad k \rightarrow \infty,
}
which implies the result.  Finally, we note that claim $(d)$ follows directly from $(b)$ and $(c)$. 
}

This lemma confirms that $G_{k,d}(1)$ gives a better bound asymptotically as $k \rightarrow \infty$ than $H_k(1)$. Illustration of this asymptotic behaviour is given in the right panel of Figure \ref{f:delta_bound_dD}. More importantly, this lemma shows the overall advantage of sampling derivatives, i.e.\ we have the following:

\cor{\label{cor:large_k_HD}
For large $k$, the set $\{ \sqrt{\mu_{n,\alpha}} D^{\alpha} \Phi_{\Omega}(\cdot -x_n) : n \in I, | \alpha |_1 \leq k \}$  forms a frame for $\rB(\Omega)$ with frame bounds satisfying \R{frame_bounds_thm}, provided 
\bes{
\delta < \frac{1}{\E} \frac{k+1}{m_{\Omega} b}.
}
}

Hence, for all dimensions $d$, the maximum allowed density $\delta$ increases linearly with the number of derivatives $k$.  Unfortunately the constant of proportionality $1/\E \approx 0.3679$ is rather small.  Indeed, it is much smaller than in the case of equispaced samples, where the corresponding constant is $\pi/2 \approx 1.5708$.  To ameliorate this gap, we will first prove an improved estimate in \S \ref{ss:bound_1D} for the case $d=1$.  Second, in \S \ref{s:perturbation} we will prove a perturbation result for nonuniform derivative sampling.

\rem{\label{r:delta_lit}
For the case $k=0$, Beurling established the sharp, sufficient condition $\delta < \pi/2$ when $\Omega$ is the unit Euclidean ball and $\abs{\cdot}_* = \abs{\cdot}_2$, provided the sampling points $\{ x_n \}_{n \in I}$ are separated.  In \textnormal{\cite{BenedettoSpiral}} and \textnormal{\cite{OlevskiiUlanovskii}} this was extended to any compact, convex and symmetric domain $\Omega$, where $\abs{\cdot}_*$ is the norm induced by the polar set of $\Omega$.  The separation condition was removed in \textnormal{\cite{AGH2DNUGS}} by incorporating weights.  To the best of our knowledge, it is an open problem to see if similar sharp results can be proved for the case of sampling with derivatives.
}

\subsubsection{Proof of Theorem \ref{t:main_frame_dD}}\label{sss:main_thm_proof}

The proof of this theorem uses the techniques of \cite{GrochenigIrregular,GrochenigIrregularExpType,GrochenigModernSamplingBook}, and more recently \cite{AGH2DNUGS}, which were applied to the $d \geq 1$ and $k=0$ case, as well as the approach in \cite{RazafinjatovoDerivs} for the $d=1$ and $k \geq 0$ case.  We first require the following three lemmas.  In what follows, we denote Euclidean ball of radius $r$ centered at $v$ by $\cB(v,r)$, and when centre does not matter, we write $\cB_r$.

\lem{
\label{l:k0_B}
Let $\mu_n = \mu_{n,0}$, where $\mu_{n,0} = \mathrm{meas}(V_n)$ is the Lebesgue measure of $V_n$.  Then
\bes{
\sum_{n \in I} \mu_n | f(x_n) |^2 \leq \exp(2 b \delta r) \| f \|^2,\quad \forall f \in \rB(\Omega),
}
where $b$ is as in \textnormal{\R{ab_star_def}} and $r>0$ is radius of the smallest ball (with arbitrary centre) with $\Omega\subseteq\cB_r$.
}
\prf{ 
Let $\cB(\omega_t,r)$ be the minimal ball such that $\Omega\subseteq\cB(\omega_t,r)$ and note that we can use the following shifting argument. For every $f \in \rB(\Omega)$, if $F\in\rB(\Omega-\omega_t)$ is defined so that $\hat{F}(\omega)=\hat{f}(\omega+\omega_t)$, then we have $|F|=|f|$ and also $\|F\|=\|f\|$. Therefore, without loss of generality, we may assume that $\Omega\subseteq\cB(0,r)$.  Since a bandlimited function is analytic, by Taylor's theorem we have
\bes{
f(x_n) = \sum_{\alpha \in \bbN^d_0} \frac{(x_n-x)^{\alpha}}{\alpha!} D^{\alpha} f(x),
}
for any $n\in I$ and $x\in\bbR^d$.  Let $c >0$ be a constant.  By the Cauchy--Schwarz inequality
\bes{
| f(x_n) |^2 \leq \sum_{\alpha \in \bbN^d_0} \frac{|(x-x_n)^{2\alpha}| c^{|\alpha|_1}}{\alpha!} \sum_{\alpha \in \bbN^d_0} \frac{c^{-|\alpha|_1}}{\alpha!} | D^{\alpha} f(x)|^2.
}
Note that $|(x-x_n)^{2\alpha}|=\prod_{j=1}^d |(x-x_n)_j|^{2\alpha_j}$, by using the multi-index notation introduced in Section \ref{s:pre}. 
Next, by the multinomial formula \R{multinomial},
\eas{
\sum_{\alpha \in \bbN^d_0} \frac{|(x-x_n)^{2\alpha}| c^{|\alpha|_1}}{\alpha!} = \sum^{\infty}_{k=0} \frac{c^k}{k!} \sum_{| \alpha |_1 = k} \frac{k!}{\alpha!} |(x-x_n)^{2\alpha}| &=  \sum^{\infty}_{k=0} \frac{c^k}{k!} | x - x_n |^{2k}_2 \\
& = \exp(c | x - x_n |^2_2).
}
By \R{ab_star_def} and the definition of $\delta$, we have $| x- x_n |_2 \leq b | x - x_n |_* \leq b \delta$.  Hence we find that
\bes{
| f(x_n) |^2 \leq \exp(c b^2 \delta^2) \sum_{\alpha \in \bbN^d_0} \frac{c^{-|\alpha|_1}}{\alpha!} | D^{\alpha} f(x)|^2.
}
Using definition of $\mu_n$ and the fact that Voronoi cells form a partition of $\bbR^d$, this now gives
\eas{
\sum_{n \in I} \mu_n | f(x_n) |^2 &\leq \exp(c b^2 \delta^2)\sum_{\alpha \in \bbN^d_0} \frac{c^{-|\alpha|_1}}{\alpha!} \sum_{n \in I} \int_{V_n} | D^{\alpha} f(x)|^2 \D x \\
&= \exp(c b^2 \delta^2)\sum_{\alpha \in \bbN^d_0} \frac{c^{-|\alpha|_1}}{\alpha!} \| D^{\alpha} f \|^2.
}
Consider the sum.  We have, from the standard properties of the Fourier transform,
\eas{
\sum_{\alpha \in \bbN^d_0} \frac{c^{-|\alpha|_1}}{\alpha!} \| D^{\alpha} f \|^2 &= \sum_{\alpha \in \bbN^d_0} \frac{c^{-|\alpha|_1}}{(2 \pi)^d \alpha!} \| \widehat{D^{\alpha} f} \|^2
\\
& = \sum_{\alpha \in \bbN^d_0} \frac{c^{-|\alpha|_1}}{(2\pi)^d \alpha!} \int_{\Omega} |\omega^{2\alpha} | | \hat{f}(\omega) |^2 \D \omega
\\
& = \frac{1}{(2\pi)^d} \int_{\Omega} \exp(| \omega |^2_2 / c ) | \hat{f}(\omega) |^2 \D \omega
\\
& \leq \exp(r^2 / c) \| f \|^2.
}
Therefore we deduce that
\bes{
\sum_{n \in I} \mu_n | f(x_n) |^2 \leq \exp(c b^2 \delta^2 + r^2 / c) \| f \|^2,
}
and setting $c = r / (b \delta)$ gives the result.
}

Later we will see that application of this lemma leads (after some additional work) to the bound $\delta < H_k(1) / (m_{\Omega} b)$.  As discussed, this does not give the best scaling as $k \rightarrow \infty$, which can be traced to the exponential growth in $\delta$ of the bound obtained in this lemma.  In order to mitigate this growth, and therefore eventually get a better density bound, we need the following result.

\lem{
\label{l:k0_B_by_splitting_domain}
Let $\mu_n = \mu_{n,0}$, where $\mu_{n,0} = \mathrm{meas}(V_n)$ is the Lebesgue measure of $V_n$.  Then
\bes{
\sum_{n \in I} \mu_n | f(x_n) |^2 \leq (1+2\sigma_d^*(b\delta m_\Omega))^d \exp(2 b \delta m_\Omega/\sigma_d^*(b\delta m_\Omega)) \| f \|^2,\quad \forall f \in \rB(\Omega),
}
where $\sigma_d^*$, $b$ and $m_{\Omega}$ are as in \textnormal{\R{sigma_star}}, \textnormal{\R{ab_star_def}} and \textnormal{\R{m_Omega_def}} respectively. 
}

\prf{ 
Let $\sigma>0$ be fixed and let us cover $\Omega$ by $R=R(\Omega,\cB_{m_\Omega/\sigma})$ Euclidean balls of radius $m_\Omega/\sigma$. By using a classical result on covering numbers, see for example \cite{FoucartRauhutCompressBook}, we have
\be{\label{covering_no}
R \leq R(\cB_{m_\Omega},\cB_{m_\Omega/\sigma}) = R(\cB_1,\cB_{1/\sigma}) \leq (1+2\sigma)^d. 
}
Let $\left\{\cB^1_{m_\Omega/\sigma},\ldots,\cB^R_{m_\Omega/\sigma}\right\}$ be the prescribed cover of $R$ balls for $\Omega$.  Using this cover, we form a partition of $\Omega$ as follows. Set $\Omega_1 = \cB^1_{m_\Omega/\sigma} \cap \Omega $, and given $\Omega_{1},\ldots,\Omega_r$, define
\bes{
\Omega_{r+1} = \left( \cB^{r+1}_{m_\Omega/\sigma} \cap \Omega \right) \backslash \bigcup^{r}_{j=1} \Omega_j.
}
This gives at most $R$ nonempty sets $\Omega_1,\ldots,\Omega_R$ which make a disjoint cover of $\Omega$. By construction, for each $j$, $\Omega_j\subseteq\cB^j_{m_\Omega/\sigma}$. Due to Lemma \ref{l:k0_B}, we know that
\be{\label{eq:lemma_use}
\sum_{n \in I} \mu_n | g(x_n) |^2 \leq \exp(2 b \delta m_\Omega/\sigma) \| g \|^2,\quad \forall g \in \rB(\Omega_j), \quad j=1,\ldots,R.
}
Since $\Omega_1,\ldots,\Omega_R$ are disjoint and $\bigcup_{j=1}^R\Omega_j=\Omega$, for each $f\in\rB(\Omega)$ we have that $\hat{f}=\sum_{j=1}^{R}\hat{f}_j$, $f=\sum_{j=1}^R f_j$ and $\|f\|^2 = \sum_{j=1}^R \|f_j\|^2$, 
where $f_j\in\rB(\Omega_j)$. Therefore we get
\bes{
\sum_{n \in I} \mu_n | f(x_n) |^2  \leq R \sum_{j=1}^{R}  \sum_{n \in I} \mu_n   | f_j(x_n) |^2 \leq (1+2\sigma)^d \exp(2 b \delta m_\Omega/\sigma) \| f \|^2.
}
Now, if we minimize the right-hand side over $\sigma>0$, we obtain
\eas{
&\sum_{n \in I} \mu_n | f(x_n) |^2  \\
&\leq \left(1+2\frac{b\delta m_\Omega+\sqrt{b\delta m_\Omega(d+b\delta m_\Omega)}}{d}\right)^d \exp\left(\frac{2 b \delta m_\Omega d}{b\delta m_\Omega+\sqrt{b\delta m_\Omega(d+b\delta m_\Omega)}} \right) \| f \|^2
}
and the result follows.
}

\rem{
Note that Lemma $\ref{l:k0_B_by_splitting_domain}$ also holds when $m_{\Omega}$ is replaced by $r$, where $r$ is defined as in Lemma $\ref{l:k0_B}$. This holds due to the same shifting argument as used in the proof of Lemma $\ref{l:k0_B}$. While this shifting argument proved to be crucial when applying Lemma $\ref{l:k0_B}$ for derivation of $\R{eq:lemma_use}$, in the results that follow such shifting argument is not of much use. This is so because in general $|D^{\alpha}F|\neq|D^{\alpha}f|$, when $F$ is defined as in the proof of Lemma $\ref{l:k0_B}$.}

\lem{
\label{l:f_Taylor_d}
For any $f\in\rB(\Omega)$, we have
\bes{
\nm{f - \sum_{n \in I} \sum_{| \alpha |_1 \leq k} \frac{1}{\alpha!} D^{\alpha}f(x_n) (\cdot-x_n)^{\alpha} \bbI_{V_n} } \leq  \min\left\{ h_k( b \delta m_\Omega), g_{k,d}( b \delta m_\Omega) \right\} \|f\|
}
where $h_k$ and $g_{k,d}$ are as in \textnormal{\R{h_k}} and \textnormal{\R{g_k,d}}, and $b$ and $m_\Omega$   are as in \textnormal{\R{ab_star_def}} and \textnormal{\R{m_Omega_def}} respectively.
}

\prf{
For $f\in\rB(\Omega)$ let $g(x)=\sum_{n \in I} \sum_{| \alpha |_1 \leq k} \frac{1}{\alpha!} D^{\alpha}f(x_n) (x-x_n)^{\alpha} \bbI_{V_n}(x)$.
Since Voronoi cells form a partition of $\bbR^d$, we have
\bes{
\nm{f - g }^2 = \sum_{n \in I} \int_{V_n} \left | f(x) - \sum_{| \alpha |_1 \leq k} \frac{1}{\alpha!} D^{\alpha}f(x_n) (x-x_n)^{\alpha} \right |^2 \D x.
}
Let $x \in V_n$.  By Taylor's theorem and the Cauchy--Schwarz inequality
\eas{
&\left | f(x) - \sum_{| \alpha |_1 \leq k} \frac{1}{\alpha!} D^{\alpha}f(x_n) (x-x_n)^{\alpha} \right |^2 \\
&=\left |  \sum_{| \alpha |_1 > k}  \frac{1}{\alpha!} D^{\alpha}f(x_n) (x-x_n)^{\alpha} \right |^2
\\
& \leq \sum_{|\alpha|_1 > k} \frac{c^{|\alpha|_1} | (x-x_n)^{2\alpha} | }{\alpha!} \sum_{|\alpha|_1 > k} \frac{c^{-|\alpha|_1}}{\alpha!} | D^{\alpha}f(x_n) |^2.
}
Note that
\eas{
\sum_{|\alpha|_1 > k} \frac{c^{|\alpha|_1} | (x-x_n)^{2\alpha} | }{\alpha!} = \sum_{r > k} \frac{c^r}{r!} | x - x_n |^{2r}_2 \leq R_k(c b^2 \delta^2),
}
where $R_k$ is as in \R{R_k}.  Hence, by Lemma \ref{l:k0_B} applied to $D^{\alpha}f\in\rB(\Omega)$, 
\eas{
\nm{f - g }^2 &\leq R_k(c b^2 \delta^2) \sum_{| \alpha |_1 > k} \frac{c^{-|\alpha|_1}}{\alpha!} \sum_{n \in I} \mu_{n} | D^{\alpha}f(x_n) |^2
\\
& \leq R_k(c b^2 \delta^2) \exp(2 b \delta m_\Omega) \sum_{| \alpha |_1 > k} \frac{c^{-|\alpha|_1}}{\alpha!} \| D^{\alpha} f \|^2.
}
Noting that 
\bes{
\sum_{| \alpha |_1 > k} \frac{c^{-|\alpha|_1}}{\alpha!} \| D^{\alpha} f \|^2 =\frac{1}{(2\pi)^d} \int_{\Omega} \sum_{|\alpha|_1 > k} \frac{c^{-|\alpha|_1}}{\alpha!} | \omega^{2 \alpha} | | \hat{f}(\omega) |^2 \D \omega \leq R_k(m_\Omega^2 / c) \| f \|^2
}
and setting $c = m_\Omega / (b \delta)$ gives
$\nm{f - g } \leq R_k(b\delta m_\Omega) \exp(b \delta m_\Omega) \| f \| = h_k(b \delta m_\Omega) \| f \|$.
Similarly, if we apply Lemma \ref{l:k0_B_by_splitting_domain} instead of Lemma \ref{l:k0_B}, we get 
\bes{
\nm{f - g } \leq (1+2\sigma_d^*(b\delta m_\Omega))^{d/2} R_k(b\delta m_\Omega) \exp(b \delta m_\Omega/\sigma_d^*(b\delta m_\Omega)) \| f \| =g_{k,d}(b \delta m_\Omega) \| f \|,
}
and the result follows.
}

Now we are ready to prove Theorem \ref{t:main_frame_dD}.

\prf{[Proof of Theorem \ref{t:main_frame_dD}]
Fix $f \in \rB(\Omega)$ and define $g (\cdot)= \sum_{n \in I} \sum_{| \alpha |_1 \leq k} \frac{1}{\alpha!} \linebreak D^{\alpha}f(x_n) (\cdot-x_n)^{\alpha} \bbI_{V_n}(\cdot)$.  For the upper bound on $\| g \|^2$ we have
\eas{
\sum_{n \in I} \int_{V_n} \left | \sum_{| \alpha |_1 \leq k} \frac{1}{\alpha!} D^{\alpha} f(x_n) (x-x_n)^{\alpha} \right |^2 \D x\leq \left ( \sum_{|\alpha|_1 \leq k} \frac{1}{\alpha!} \right ) \sum_{n \in I} \sum_{| \alpha |_1 \leq k } \mu_{n,\alpha} |D^{\alpha} f(x_n) |^2.
}
By the multinomial formula
\bes{
\sum_{|\alpha|_1 \leq k} \frac{1}{\alpha!} = \sum^{k}_{l=0} \frac{1}{l!} \sum_{|\alpha|_1 = l} \frac{l!}{\alpha!}= \sum^{k}_{l=0} \frac{d^l}{l!} \leq \E^{d}.
}
Using this we get
\ea{
\label{lower_skip}
\sum_{n \in I} \sum_{| \alpha |_1 \leq k } \mu_{n,\alpha} |D^{\alpha} f(x_n) |^2 \geq \E^{-d} \| g \|^2 \geq \E^{-d} \left ( \| f \| - \| f - g \| \right )^2.
}
Lemma \ref{l:f_Taylor_d} now gives the lower bound.  Next, we address the upper bound.   Note that
\bes{
\mu_{n,\alpha} \leq \frac{1}{\alpha!} \sup_{x \in V_n} |(x-x_n)^{2 \alpha}| \mu_{n,0}.
}
Moreover $|(x-x_n)^{2 \alpha}| \leq | x-x_n |^{2 | \alpha |_1}_{\infty} \leq | x - x_n |^{2 |\alpha|_1}_{2} \leq (b \delta)^{2 | \alpha |_1}$.
Hence, Lemma \ref{l:k0_B} gives
\eas{
\sum_{n \in I} \sum_{| \alpha |_1 \leq k} \mu_{n,\alpha} | D^{\alpha}f(x_n) |^2 \leq \exp(2 m_{\Omega} b \delta) \sum_{| \alpha |_1 \leq k} \frac{(b \delta)^{2 | \alpha |_1}}{\alpha!} \| D^{\alpha} f \|^2.
}
Arguing in the same way now yields
\bes{
\sum_{n \in I} \sum_{| \alpha |_1 \leq k} \mu_{n,\alpha} | D^{\alpha}f(x_n) |^2 \leq \exp(2 m_{\Omega} b \delta + (m_{\Omega} b \delta)^2) \| f \|^2,
}
as required.
}

\subsection{The univariate case}\label{ss:bound_1D}
In the one-dimensional setting it is possible to improve the bound derived in Theorem \ref{t:main_frame_dD} somewhat using so-called Wirtinger inequalities.  See \cite{GrochenigIrregular} for the case $k=0$ and \cite{RazafinjatovoDerivs} for $k=1$.  

Throughout this section $\Omega \subseteq \bbR$ is compact and $\{ x_n \}_{n \in \bbZ}$ is a set of sampling points in $\bbR$, indexed over $\bbZ$.  We assume the points are ordered so that $x_n < x_{n+1}$, $\forall n \in \bbZ$.  As before, we let
\be{\label{delta_1d_def}
\delta = \sup_{x \in \bbR} \inf_{n \in \bbZ} | x - x_n |,
}
where $\abs{\cdot}$ denotes the absolute value.  Note that the Voronoi cells $V_n$ are the intervals
\bes{
V_n = [z_n,z_{n+1}],\quad z_n = \frac{x_n+x_{n-1}}{2},\qquad n \in \bbZ.
}
As stated above, we shall use Wirtinger inequalities to derive bounds for $\delta$.  Specifically, for $k \in \bbN$, let $c_k > 0$ be the minimal constant such that
\be{
\label{Wirtinger_higher}
\int^{b}_{a} | f(x) |^2 \D x \leq (c_k)^{2k} (b-a)^{2k} \int^{b}_{a} | f^{(k)}(x) |^2 \D x,
}
for all $f \in \rH^k(a,b)$, the $k^{\rth}$ Sobolev space, satisfying
\bes{
f(a) = f'(a) = \ldots = f^{(k-1)}(a) = 0\ \quad \mbox{or}\ \quad f(b) = f'(b) = \ldots = f^{(k-1)}(b) = 0.
}

\thm{
\label{t:main_fram_1D}
Suppose that the weights are
\bes{
\mu_{n,l} = \frac{1}{l!} \int_{V_n} (x-x_n)^{2 l} \D x = \frac{(z_{n+1}-x_n)^{2l+1} - (z_n -x_n)^{2l+1} }{l! (2l+1)},\ l=0,\ldots,k,\ n \in \bbZ,
}
and let $\delta$ be as in \textnormal{\R{delta_1d_def}}.
If
\be{
\label{delta_1D}
\delta < \frac{C(k)}{ m_{\Omega}},\qquad C(k) = \frac{1}{c_{k+1}},
}
where $c_k$ is as in \textnormal{\R{Wirtinger_higher}}, 
then
\bes{
A \| f \|^2 \leq \sum_{n \in \bbZ} \sum^{k}_{l=0} \mu_{n,l} | f^{(l)}(x_n) |^2 \leq B \| f \|^2,\quad \forall f \in \rB(\Omega),
} 
where
\bes{
A \geq \E^{-1} \left ( 1 - (c_{k+1}\delta m_{\Omega})^{k+1} \right )^2,\quad B \leq \left ( 1 +  2\delta m_{\Omega}/\pi \right )^2 \E^{(\delta m_{\Omega} )^2 }.
}
Equivalently, the set $\{ \sqrt{\mu_{n,l}} \frac{\D^l}{\D x^l} \Phi_{\Omega}(\cdot -x_n) : n \in \bbZ, l=0,\ldots,k \}$ forms a frame for $\rB(\Omega)$ with bounds $A$ and $B$. 
}

In the following section we examine the constants $c_k$ and conclude by discussion with the improvement offered by this theorem over the multivariate result Theorem \ref{t:main_frame_dD}.

\prf{
We follow the arguments of \cite{GrochenigIrregularExpType,RazafinjatovoDerivs}.  Let $g(x) = \sum_{n \in \bbZ} \sum^{k}_{l=0} \frac{1}{l!} f^{(l)}(x_n) (x-x_n)^l \bbI_{V_n}(t)$. 
Then 
\eas{
\| f - g \|^2 &=\sum_{n \in \bbZ} \int^{z_{n+1}}_{z_n} \left | f(x) - \sum^{k}_{l=0} \frac{1}{l!} f^{(l)}(x_n) (x-x_n)^l \right |^2 \D x.
\\
& = \sum_{n \in \bbZ} \left ( \int^{z_{n+1}}_{x_n} + \int^{x_n}_{z_n} \right ) \left | f(x) - \sum^{k}_{l=0} \frac{1}{l!} f^{(l)}(x_n) (x-x_n)^l \right |^2 \D x.
}
The function $f(x) - \sum^{k}_{l=0} \frac{1}{l!} f^{(l)}(x_n) (x-x_n)^l$ vanishes, along with its first $k$ derivatives, at $x = x_n$.  Applying \R{Wirtinger_higher} to each integral and noting that $|z_{n+1} - x_n | \leq \delta$ and $ | x_n - z_n | \leq \delta$ gives
\bes{
\| f - g \|^2 \leq (c_{k+1}\delta)^{2k+2} \| f^{(k+1)} \|^2.
}
Applying Bernstein's inequality \R{Bernstein} we deduce that $\| f - g \| \leq (c_{k+1}\delta m_{\Omega})^{k+1} \| f \|$, and therefore
\bes{
\left ( 1 - (c_{k+1}\delta m_{\Omega})^{k+1} \right ) \| f \| \leq \| g \| \leq \left ( 1 + (c_{k+1}\delta m_{\Omega})^{k+1} \right ) \| f \|.
}
We now use this and \R{lower_skip} to get the estimate for $A$.  For the bound on $B$, we argue similarly to the proof of Theorem \ref{t:main_frame_dD}.  We have
\bes{
\mu_{n,l} \leq \frac{1}{l!} \sup_{x \in V_n} | x-x_n |^{2l} \mu_{n,0} \leq \frac{1}{l!} \delta^{2l} \mu_{n,0}.
}
Hence
\bes{
\sum_{n \in \bbZ} \sum^{k}_{l=0} \mu_{n,l} | f^{(l)}(x_n) |^2 \leq \sum^{k}_{l=0} \frac{\delta^{2l}}{l!} \sum_{n \in \bbZ} \mu_{n,0} | f^{(l)}(x_n) |^2.
}
Gr\"ochenig's result \cite{GrochenigIrregular} for $k=0$ and $d=1$ gives that $\sum_{n \in \bbZ} \mu_{n,0} | g(x_n) |^2 \leq (1+2\delta m_{\Omega}/\pi)^2 \| g \|^2$, $\forall g \in \rB(\Omega)$.  By this and Bernstein's inequality, we deduce that
\bes{
\sum_{n \in \bbZ} \sum^{k}_{l=0} \mu_{n,l} | f^{(l)}(x_n) |^2 \leq (1+2\delta m_{\Omega}/\pi)^2 \sum^{k}_{l=0} \frac{(\delta m_{\Omega})^{2l}}{l!}  \| f \|^2.
}
Since $\sum^{k}_{l=0} \frac{(\delta m_{\Omega})^{2l}}{l!}=\exp((\delta m_{\Omega})^2)$, the upper bound follows.
}

Observe that for $k=0$, i.e.\ the classical nonuniform sampling problem without derivatives, \R{delta_1D} reduces to $\delta < \frac{\pi }{2m_{\Omega}}$ since $c_1 = \frac{2}{\pi}$ \cite{GrochenigIrregular}.  This is in agreement with the result of Gr\"ochenig \cite{GrochenigIrregular} (note that the definition of $\delta$ used therein is precisely twice the definition we use here).  This result is sharp, and says that one must sample at a rate just above the Nyquist rate $\frac{\pi}{2 m_{\Omega}}$.

\subsubsection{The magnitude of $c_k$}\label{sss:Wirtinger}

We now consider the case $k \geq 1$.  The main issue is the magnitude of the constant $c_k$ of Wirtinger's inequality \R{Wirtinger_higher}.  We first note the following:

\lem{
Consider the polyharmonic eigenvalue problem
\be{
\label{eig_problem}
(-1)^k g^{(2k)} = \lambda g,\quad g(0)=\ldots=g^{(k-1)}(0) = g^{(k)}(1) = \ldots g^{2k-1}(1) = 0.
}
This problem has a countable basis of positive eigenvalues $0 < \lambda^{(k)}_1 < \lambda^{(k)}_2 < \ldots$.  Moreover, the best constant $c_k$ in the inequality \textnormal{\R{Wirtinger_higher}} is precisely $(\lambda^{(k)}_1)^{-\frac{1}{2k}}$.
}
\prf{
It is well known that \R{eig_problem} has a countable spectrum with eigenfunctions $\{ \phi_n \}^{\infty}_{n=1}$ forming an orthonormal basis of $L^2(0,1)$ \cite{Naimark1}.  It is straightforward to see that \R{eig_problem} has only strictly positive eigenvalues.  Now let $f \in \rH^k(0,1)$ satisfy $f(0)=\ldots=f^{(k-1)}(0)=0$.  Then
\bes{
\ip{f}{\phi_n} = \frac{(-1)^k}{\lambda^{(k)}_n} \ip{f}{\phi^{(2k)}_n} = \frac{1}{\lambda^{(k)}_n} \ip{f^{(k)}}{\phi^{(k)}_n}.
}
In particular, if $f = \phi_n$, then $\| \phi_n \|^2 = \frac{1}{\lambda^{(k)}_n} \| \phi^{(k)}_n \|^2$.  Let $\psi_n = \frac{1}{\sqrt{\lambda^{(k)}_n}} \phi^{(k)}_n$, so that $\| \psi_n \| =1$.  The set $\{ \psi_n \}^{\infty}_{n=1}$ is precisely the set of eigenfunctions of the problem
\bes{
(-1)^k g^{(2k)} = \lambda g,\quad g^{(k)}(0)=\ldots=g^{(2k-1)}(0) = g(1) = \ldots g^{k-1}(1) = 0.
}
In particular, they form an orthonormal basis of $\rL^2(0,1)$.  Since $\ip{f}{\phi_n} = \frac{1}{\sqrt{\lambda^{(k)}_n}} \ip{f^{(k)}}{\psi_n}$,
it follows from Parseval's identity that
\bes{
\| f \|^2 = \sum_{n} | \ip{f}{\phi_n} |^2 = \sum_{n} \frac{1}{\lambda^{(k)}_n} | \ip{f^{(k)}}{\psi_n} |^2 \leq \frac{1}{\lambda^{(k)}_1} \sum_{n} | \ip{f^{(k)}}{\psi_n} |^2 = \frac{1}{\lambda^{(k)}_1} \| f^{(k)} \|^2,
}
by completeness.  Thus $\| f \|^2 \leq 1/\lambda^{(k)}_1 \| f^{(k)} \|^2$, and this bound is sharp since we may set $f = \phi_1$.  By a change of variables, we get that $(c_k)^{2k} = 1/\lambda^{(k)}_1$, as required.
}

This means we can determine the constant $c_k$ by finding the eigenvalues of \R{eig_problem}.  When $k=1$, the eigenvalues of \R{eig_problem} are $(\pi/2 + n \pi )^2$, $n \in \bbN_0$.  Hence $\lambda^{(1)}_1 = \pi^2/4$ and $c_1 = 2 / \pi$, as stated.  Unfortunately, for $k \geq 2$ no explicit expression exists for the eigenvalues, so we resort to numerical computation.  For $k \geq 2$, write $\lambda = \tau^{2k}$ for $\tau > 0$.  The general solution of \R{eig_problem} can be written as
\bes{
g(x) = \sum^{2k-1}_{s=0} b_s \E^{\I z^s \tau x},
}
where $z = \E^{\I \pi / k}$ and $b_s \in \bbC$ are coefficients.  Enforcing the boundary conditions results in a linear system of equations
\eas{
\sum^{2k-1}_{s=0} (\I z^s \tau)^r b_s = 0,\quad \sum^{2k-1}_{s=0} (\I z^s \tau)^{k+r} \E^{\I z^s \tau} b_s = 0, \qquad r=0,\ldots,k-1.
}
In matrix form, we have $A(\tau) b = 0$, where $A(\tau) \in \bbC^{2k \times 2k}, b = (b_0,\ldots,b_{2k-1})^{\top}$.  Hence the minimal eigenvalue $\lambda^{(k)}_1 = (\tau^{(k)}_1)^{2k}$, and therefore $c_k = 1 / \tau^{(k)}_1$, where $\tau^{(k)}_1$ is the first positive root of the function $D(\tau) = \det(A(\tau))$.
In the case $k=2$, we have $D(\tau) = 8 \I \tau^6 \left (  1+ \cos(\tau) \cosh(\tau) \right )$,
and numerical computation finds that $\tau^{(2)}_1 = 1.8751$ (see also \cite{RazafinjatovoDerivs}).

\begin{table}
\begin{center}
\scalebox{0.85}{
\begin{tabular}{|c|c|c|c|c|c|c|c|c|c|c|c|c|c|c|c|}
\hline
$k$ & 1 & 2 & 3 & 4 & 5 & 6 & 7 & 8 & 9 & 10  \\ \hline
$c_k$    &  0.6366 & 0.5333 & 0.4495 & 0.3861 &  0.3376 & 0.2997 & 0.2694 & 0.2446 & 0.2240 & 0.2066   \\ \hline
$1/c_k$ &  1.5708 & 1.8751 & 2.2248 & 2.5903 & 2.9621 & 3.3367 & 3.7125 & 4.0888 & 4.4652 & 4.8415   \\ \hline
\end{tabular}
}
\end{center}\caption{The values $c_k$ and $1/c_k$ for $k=1,2,\ldots,10$.  These values were calculated in high precision using \textit{Mathematica}.}
\label{tab:WirtingerConst}
\end{table}

In Table \ref{tab:WirtingerConst} we compute $\tau^{(k)}_1 = 1/c_k$ and $c_k$ for $k=1,\ldots,10$.  As is evident the values $1/c_k$, grow approximately linearly in $k$ for large $k$.  Linear regression on the computed values gives that $1/c_k \approx 1.1458 + 0.3674 k$ for large $k$.  Note that $1/\E = 0.3679$.  We therefore conjecture that
\be{
\label{Bottcher_Conj}
\frac{1}{c_k} \sim \frac{k}{\E},\quad k \rightarrow \infty.
}
We remark in passing that the large $k$ asymptotics for the optimal constant in a variant of Wirtinger's inequality where $f$ and its derivatives vanish at both endpoints has been derived by B\"ottcher \& Widom \cite{BottcherWidomWirtinger}.  We expect a similar approach can be applied to \R{Wirtinger_higher} to obtain \R{Bottcher_Conj}.

We can now compare Theorem \ref{t:main_fram_1D} with the multivariate result Theorem \ref{t:main_frame_dD}.  In Table \ref{tab:delta_bound} we give the numerical values for the constant $C(k)$ arising from both theorems, where $\delta < C(k) / m_{\Omega}$ is the required condition on $\delta$.  The univariate bound is evidently superior for all values of $k$ considered.   However, the bounds behave the same asymptotically, since both Theorem \ref{t:main_frame_dD} and Theorem \ref{t:main_fram_1D} give $C(k) \sim 1/\E (k+1) \approx 0.3679 (k+1)$ for large $k$ (recall Corollary \ref{cor:large_k_HD}).  In Table \ref{tab:delta_bound} we also compare Theorem \ref{t:main_fram_1D} to the bound derived in \cite[Thm.\ 1]{RazafinjatovoDerivs} (note that the value $1.8751$ for $k=1$ was also provided in \cite{RazafinjatovoDerivs} using Wirtinger's inequality arguments as we do above).  Unfortunately, the improvement obtained from Theorem \ref{t:main_fram_1D} is only marginal.  In particular, both bounds are asymptotic to $1/\E (k+1)$ for large $k$, and therefore (we expect) a long way from being sharp (recall that the condition for equispaced samples is $\delta \leq \pi /2 (k+1)$).
We conclude that although Wirtinger's inequality obtains a sharp bound for $k=0$, it is of little use in getting superior bounds for $k \geq 1$.

\begin{table}
\begin{center}
\scalebox{0.86}{
\begin{tabular}{|c|c|c|c|c|c|c|c|c|c|c|c|c|c|c|c|}
\hline
$k$ & 0 & 1 & 2 & 3 & 4 & 5 & 6 & 7 & 8 & 9  \\ \hline
(a)  &  0.4812 & 0.8141 & 1.1268 & 1.4304 & 1.7890 & 2.1535 & 2.5186 & 2.8842 & 3.2501 & 3.6163   \\ \hline
 (b)      &  1.5708 & 1.8751 & 2.2248 & 2.5903 & 2.9621 & 3.3367 & 3.7125 & 4.0888 & 4.4652 & 4.8415   \\ \hline
(c) & 1.4142 & 1.8612 & 2.2209 & 2.5886 & 2.9612 & 3.3361 & 3.7121 & 4.0885 &
4.4650 & 4.8413 \\ \hline
\end{tabular}
}
\end{center}\caption{The constant $C(k)$ obtained from \textnormal{(a)} Theorem \textnormal{\ref{t:main_frame_dD}} for the case $d=1$, \textnormal{(b)} Theorem \textnormal{\ref{t:main_fram_1D}} and \textnormal{(c)} \textnormal{\cite[Thm.\ 1]{RazafinjatovoDerivs}}. }
\label{tab:delta_bound}
\end{table}

\subsection{Line-by-line sampling}\label{ss:bound_tensor}

In some applications, not least seismology, the unknown function $f$ depends on a spatial variable $z \in \bbR^{d-1}$ and a temporal variable $t \in \bbR$.  Sensors are placed at fixed locations $\{ z_n \}_{n \in I} \subseteq \bbR^{d-1}$, where $d=2,3$, in physical space, and measurements are taken at times $\{ t_{m,n} \}_{m \in J}$.  In particular, different sensors may take measurements at different times.  This gives the set of samples
\bes{
D^{\alpha}_z f(z_n,t_{m,n}),\quad  n \in I, m \in J, | \alpha |_1 \leq k.
}
Note that $D^{\alpha}_z = \partial^{\alpha_1}_{z_1} \cdots \partial^{\alpha_{d-1}}_{z_{d-1}}$ is the partial derivative with respect to $z$ only.  We do not measure any temporal derivatives.

Let $x = (z,t) \in \bbR^d$ and write $f(z,t) = f(x)$.  We shall assume that $f \in \rB(\Omega)$ and moreover that $\Omega = \Omega_z \times \Omega_t$ for  $\Omega_z \subseteq \hat{\bbR}^{d-1}$ and $\Omega_t \subseteq \hat{\bbR}$.  Let
\bes{
\delta_z = \sup_{z \in \bbR^{d-1}} \inf_{n \in I} | z-z_n |_*,\qquad \delta_{t} = \sup_{n \in I} \sup_{t \in \bbR} \inf_{m \in J} | t - t_{m,n} |,
}
and write $V_n \subseteq \bbR^{d-1}$ for the Voronoi cells of the sampling points $\{ z_n \}_{n \in I}$ with respect to the $\abs{\cdot}_*$ norm.  We now have the following result.  Note that this is a straightforward extension of a result of Strohmer \cite{StrohmerLevinson} (see also \cite{GrochenigModernSamplingBook}) to the case of derivatives and $d \geq 3$.

\prop{
\label{p:main_tensor_1}
Suppose that the weights 
\bes{
\mu_{m,n,\alpha} = \frac{ t_{m+1,n} - t_{m,n}}{2 \alpha!} \int_{V_n} (z-z_n)^{2 \alpha} \D z .
}
If
\bes{
\delta_t < \frac{2}{\pi m_{\Omega_t}},\quad \delta_z < \frac{C(k,d)}{m_{\Omega_z} b},\qquad C(k,d) =  \left \{ \begin{array}{ll} 1/c_{k+1} & d=2 \\ \max\left\{H_k(1),G_{k,d}(1)\right\} & d \geq 3  \end{array} \right .,
}
then
\eas{
\left (1-\frac{2 \delta_t m_{\Omega_t}}{\pi} \right )^2 A_z \| f \|^2 & \leq \sum_{m \in J} \sum_{n \in I} \sum_{| \alpha |_1 \leq k} \mu_{m,n,\alpha} | D^{\alpha}_z f(z_n,t_{m,n}) |^2 \\
&\leq \left (1+\frac{2 \delta_t m_{\Omega_t}}{\pi} \right )^2 B_z \| f \|^2,
} 
for all $f \in \rB(\Omega)$, where $A_z$ and $B_z$ satisfy
\bes{
A_z \geq  \E^{-1} \left ( 1 - (m_{\Omega_z} c_{k+1}\delta_z )^{k+1} \right )^2,\quad B_z \leq \left ( 1 + 2 m_{\Omega_z} \delta_z  /\pi \right )^2 \E^{(m_{\Omega_z} \delta_z )^2 },\quad d = 2,
}
with $c_k$ as in \textnormal{\R{Wirtinger_higher}}, or
\bes{
A_z \geq \E^{-d} \left ( 1 - \min\left\{h_k(m_{\Omega_z} b \delta_z), g_{k,d}(m_{\Omega_z} b \delta_z) \right\}\right )^2,\quad B_z \leq\exp(2 m_{\Omega_z} b \delta_z+(m_{\Omega_z} b \delta_z)^2),
}
for $d\geq3$, with $h_k$ and $g_{k,d}$ as in \textnormal{\R{h_k}} and \textnormal{\R{g_k,d}} with inverse functions $H_k$ and $G_{k,d}$ respectively. 
Equivalently, the set $\{ \sqrt{\mu_{m,n,\alpha}} D^{\alpha}_z \Phi_{\Omega}(\cdot -x_{n,m}) : n \in I, m \in J, | \alpha |_1 \leq k \}$, where $x_{n,m} = (z_n,t_m)$, forms a frame for $\rB(\Omega)$ with bounds
\bes{
A \geq \left (1-2 \delta_t m_{\Omega_t}/ \pi \right )^2A_z ,\qquad B \leq \left (1+2 \delta_t m_{\Omega_t}/ \pi \right )^2 B_z.
}
}
\prf{
Gr\"ochenig's original, one-dimensional, derivative-free result \cite{GrochenigIrregular} gives that
\eas{
\left (1-2 \delta_t m_{\Omega_t}/ \pi \right )^2 \int_{\bbR} | f(z,t) |^2 \D t &\leq \sum_{m \in J} \frac{t_{m+1,n}-t_{m,n}}{2} |  f(z,t_{m,n}) |^2 \\
&\leq \left (1+2 \delta_t m_{\Omega_t}/ \pi \right )^2 \int_{\bbR} | f(z,t) |^2 \D t.
}
Hence, if $g(z) = \sqrt{\int_{\bbR} | f(z,t) |^2 \D t }$ and $\tilde{\mu}_{n,\alpha} = \frac{1}{ \alpha!} \int_{V_n} (z-z_n)^{2 \alpha} \D z $ then
\eas{
\left(1- \frac{2 \delta_t m_{\Omega_t}}{\pi} \right)^2  \sum_{n \in I} \sum_{| \alpha |_1 \leq k} \tilde{\mu}_{n,\alpha} | D^{\alpha}_z g(z_n) |^2 &\leq \sum_{m \in J} \sum_{n \in I} \sum_{| \alpha |_1 \leq k} \mu_{m,n,\alpha} | D^{\alpha}_z f(z_n,t_{m,n}) |^2
\\
& \leq \! \left( \! 1 \! + \! \frac{2 \delta_t m_{\Omega_t}}{\pi} \! \right)^2 \! \sum_{n \in I} \sum_{| \alpha |_1 \leq k} \tilde{\mu}_{n,\alpha} | D^{\alpha}_z g(z_n) |^2
}
and, to get the result, we now apply Theorem \ref{t:main_frame_dD} ($d\geq 3$) or Theorem \ref{t:main_fram_1D} ($d=2$) to the sum and note that $\int_{\bbR^{d-1}} | g(z) |^2 \D z = \| f \|^2$. 
}

This proposition implies the following.  With the above type of scheme, for stable sampling one requires (i) the usual derivative-free density for univariate nonuniform sampling in the time variable, i.e.\ $\delta_t < \frac{2}{\pi m_{\Omega_t}}$, and (ii) a density in the space variable depending on the number of derivatives.

\subsection{A multivariate perturbation result with derivatives}\label{s:perturbation}
The results proved thus far give explicit guarantees for nonuniform derivatives sampling.  However, the conditions on the density $\delta$ are more stringent than those required for uniform samples.  We now show that nonuniform sampling is possible with larger gaps under appropriate conditions.

\thm{
\label{t:perturb}
Suppose that $\{ x_n \}_{n \in I} \subseteq \bbR^d$ and $\mu_{n,\alpha} >0$, $n \in I$, $| \alpha |_1 \leq k$, are such that \R{samples_sum_stab} holds with constants $A,B>0$.  Let $\{ \tilde{x}_n \}_{n \in I} \subseteq \bbR^d$ be such that 
\be{
\label{epsilon_def}
\epsilon = \sup_{n \in I} | \tilde{x}_n - x_n |_{*} < \frac{ \log ( 1+\sqrt{A/B} )}{m_{\Omega} b},
}
then
\bes{
\tilde{A} \| f \|^2 \leq \sum_{n \in I} \sum_{| \alpha |_1 \leq k} \mu_{n,\alpha} | D^{\alpha} f(\tilde{x}_n) |^2 \leq \tilde{B} \| f \|^2,
}
where
\bes{
\tilde{A} \geq \left ( \sqrt{A} - \sqrt{B} \left ( \exp ( m_{\Omega} b \epsilon ) - 1 \right ) \right )^2,\qquad \tilde{B} \leq B \exp ( 2 m_{\Omega} b \epsilon ).
}
That is, if the set $\{ \sqrt{\mu_{n,\alpha}} D^{\alpha} \Phi_{\Omega}(\cdot -x_n) : n \in I, | \alpha |_1 \leq k \}$ forms a frame  for $\rB(\Omega)$ with bounds $A$ and $B$, then the set $\{ \sqrt{\mu_{n,\alpha}} D^{\alpha} \Phi_{\Omega}(\cdot -\tilde{x}_n) : n \in I, | \alpha |_1 \leq k \}$ forms a frame with bounds $\tilde{A}$ and $\tilde{B}$.
}
\prf{
The proof is similar to those of the earlier results.  Note first that by Minkowski inequality
\eas{
&\sqrt{\sum_{n \in I} \sum_{|\alpha|_1 \leq k} \mu_{n,\alpha} | D^{\alpha} f(\tilde{x}_n) |^2 } \\
&\geq  \sqrt{\sum_{n \in I} \sum_{|\alpha|_1 \leq k} \mu_{n,\alpha} | D^{\alpha} f(x_n) |^2 } - \sqrt{\sum_{n \in I} \sum_{|\alpha|_1 \leq k} \mu_{n,\alpha} | D^{\alpha} f(x_n)-D^{\alpha} f( \tilde{x}_n) |^2 } .
}
By identical arguments to those used in \S \ref{ss:multivariate_thms}, we have
\eas{
| g(x_n) - g(\tilde{x}_n) |^2 & \leq \left ( \exp(cb^2 \epsilon^2)-1 \right ) \sum_{|\beta|_1 > 0} \frac{c^{-|\beta|_1}}{\beta!} | D^{\beta} g(x_n) |^2,
}
for any function $g \in \rB(\Omega)$.  Using this, we deduce that
\eas{
&\sum_{n \in I} \sum_{|\alpha|_1 \leq k} \mu_{n,\alpha} | D^{\alpha} f(x_n)-D^{\alpha} f( \tilde{x}_n) |^2 \\
&\leq \left ( \exp(cb^2 \epsilon^2)-1 \right ) \sum_{| \beta |_1 > 0} \frac{c^{-|\beta|_1}}{\beta!} \sum_{n \in I} \sum_{|\alpha|_1 \leq k} \mu_{n,\alpha} | D^{\alpha} D^{\beta} f(x_n) |^2
\\
&\leq B \left ( \exp(cb^2 \epsilon^2)-1 \right ) \sum_{|\beta|_1 > 0} \frac{c^{-|\beta|_1}}{\beta!} \| D^{\beta} f \|^2
\\
 & \leq B \left ( \exp(cb^2 \epsilon^2)-1 \right )  \left ( \exp((m_\Omega)^2 / c) -1 \right ) \| f \|^2.
}
Setting $c = m_{\Omega} / (b \epsilon)$ gives $\tilde{A} \geq \left ( \sqrt{A} - \sqrt{B} \left ( \exp ( m_{\Omega} b \epsilon ) - 1 \right ) \right )^2$.  Hence, $\tilde{A} > 0$ provided that $\sqrt{A} - \sqrt{B} \left ( \exp ( m_{\Omega} b \epsilon ) - 1 \right ) > 0$.
Now, rearranging gives \R{epsilon_def}.  The upper bound for $\tilde{B}$ follows similarly.
}

As with the previous results, the right-hand side \R{epsilon_def} is dimensionless whenever $\Omega$ is contained in the unit ball and $\abs{\cdot}_*=\abs{\cdot}_q$, $1\leq q\leq 2$.  Now suppose for simplicity that $\Omega \subseteq [-1,1]^d$.  Then the points $x_n = (k+1) n \pi$, $n \in \bbZ^d$, give rise to a stable set of sampling (this is due to the fact that they give rise to a Riesz basis for $\Omega = [-1,1]^d$, and therefore a frame when $\Omega \subseteq [-1,1]^d$).  This theorem therefore allows for nonuniform samples with gaps roughly on the size of $k$, provided the sampling points $\tilde{x}_n$ are within $\epsilon$ of the $x_n$.  An issue with this result is that the ratio $A/B$ is liable to decrease with both $k$ and $d$. Hence, the maximal allowed $\epsilon$ may be rather small in practice.

In \cite[Cor.\ 6.1]{BaileySampRecov}, a multivariate perturbation result for the case $k=0$ with $x_n = n \pi$ was derived based on similar arguments.  In our notation, the result proved therein corresponds to the case $p = q = \infty$.  The precise condition given is $\epsilon < \log(2) / d$, which is equivalent to \R{epsilon_def} with $k=0$.  Note that Sun \& Zhou \cite{SunZhouMultivarTrigPerturb} also prove a perturbation result in the same setting $p = q = \infty$, but based on expanding in Laplace--Neumann eigenfunctions, rather than Taylor series (this is similar to the proof of the original Kadec-1/4 theorem).  Their constant is somewhat smaller than $\log(2) / d$ for finite $d$, but, as discussed in \cite{BaileySampRecov}, it is asymptotic to $\log(2)/d$ as $d \rightarrow \infty$.  The generalizations of these results offered by Theorem \ref{t:perturb} are: (i) flexibility over the choice of domain $\Omega$---in particular, a dimension-independent bound for appropriate $\Omega$ and $\abs{\cdot}_{*}$, and (ii) the case $k \neq 0$.

In \cite{AcostaAkramKrishtal_2009}, perturbation results are proved for a more general sampling model that includes derivatives sampling of bandlimited functions as a special case. However, for this particular case \cite[Thm 3.8]{AcostaAkramKrishtal_2009}, the perturbation bound is not explicit, and additionally, it assumes separation of the sampling points.

\section{Univariate nonuniform bunched sampling}\label{s:bunched_samp}
We now consider nonuniform sampling with sampling points clustered in bunches.  Given the difficulty of polynomial interpolation for $d \geq 2$ dimensions, we consider the univariate case only.

\subsection{Setup}
Assume that we are given samples at the points $\{x_{n,0}\}_{n\in I}\subseteq\bbR$ which are $\delta$-dense  
\bes{
\delta = \sup_{x \in \bbR} \inf_{n \in I} | x - x_{n,0} |,
}
and also for each $n\in I$  we are given $s$ additional samples at the distinct points  
\be{\label{condition}
x_{n,m} \in [x_{n,0}-h_n,x_{n,0}+h_n] \subseteq V_{n},\quad m=1,\ldots,s,
} 
where $V_{n}$ is the Voronoi region associated to $x_{n,0}$. If we denote $h=\sup_{n\in I}h_n$, then $h=\tau\delta$ for some positive constant $\tau\leq1$.
Therefore, in each $h$-vicinity of $x_{n,0}$, there are $s$ additional sampling points. 
We shall call such a sampling sequence $\{x_{n,m}\}_{n\in I,0\leq m\leq s}$ a \textit{bunched sampling set} with density $\delta$ and bunch width $h$.  

Much as in the case of derivatives sampling, in bunched sampling, we expect that a larger $\delta$ is possible if there are multiple sample points around each $x_{n,0}$.  As discussed in the introduction, it is useful to have this type of sampling scheme in the situations where we must allow for bigger distances between sampling sensors due to some natural constraints.

\subsection{Bunched sampling and fusion frames}\label{ss:bunched_fusion}

In nonuniform derivative sampling, we showed the existence of a particular frame to establish stable sampling.  In the case of bunched sampling, we will first show the existence of a particular \textit{fusion} frame \cite{CasazzaKutyniokFusion,CasazzaKutyniokLi}.  We recall that a non-orthogonal fusion frame \cite{CasazzaNonortFusionFrames} for a Hilbert space $\rH$ is a set of positive scalars $\{ v_n \}_{n \in I}$ and non-orthogonal projections $\{ \cP_n \}_{n \in I}$, each with closed range, satisfying
\bes{
A \| f \|^2 \leq \sum_{n \in I} v_n \| \cP_n f \|^2 \leq B \| f \|^2,\quad \forall f \in \rH.
}
Much like a frame, the associated \textit{fusion frame operator} $\cS : \rH \rightarrow \rH$ given by
\bes{
\cS f = \sum_{n \in I} \cP^*_n \cP_n f
}
is self-adjoint and invertible.  Thus, any $f \in \rH$ can be recovered stably from the data $\{ \cP_n f \}_{n \in I}$.  In practice, if the projections have finite-dimensional ranges, the reconstruction can be carried out via generalized sampling \cite{BAMGACHNonuniform1D, BAACHShannon,BAACHOptimality}, for example.

Given the bunched sampling set  $\{x_{n,m}\}_{n\in I,0\leq m\leq s}$ and associated Voronoi regions $\{V_n\}_{n\in I}$, for each $n\in I$ we define the subspace
\bes{
W_n = \left\{ g \in \rL^2(\bbR) : \supp(g) \subseteq V_n  \right\}
}
and also for any $f\in\rB(\Omega)$ we define the operator
\be{\label{interp_operator}
\cP_n (f) = p_n(f) \bbI_{V_n}
}
where $p_n(f) \in \bbP_s$ is the unique interpolating polynomial of degree $s$ such that 
\bes{
p_n(f)(x_{n,m}) = f(x_{n,m}), \quad m=0,\ldots,s.
}
The bounded linear operator $\cP_n:\rB(\Omega)\rightarrow W_n$  is a non-orthogonal projection, i.e.\ $\cP_n^2=\cP_n$ by uniqueness of the interpolating polynomial.  Hence, if  there exist $A,B>0$ such that for all $f\in\rB(\Omega)$ 
\bes{
A\|f\|^2 \leq \sum_{n\in I} \left \| \cP_n(f) \right \|^2 \leq B \|f\|^2,
}
then $\left\{ \cP_n \right\}_{n\in I}$ is a non-orthogonal fusion frame for $\rB(\Omega)$ with weights $v_n = 1$.  Our main result gives conditions for this to be the case:

\thm{\label{t:bunched_samp}
Suppose that $\{x_{n,m}\}_{n\in I,0\leq m\leq s}\subseteq\bbR$ is a bunched sampling set with density $\delta$ and bunch width $h=\tau\delta$, where $\tau\in(0,1]$. If 
\be{\label{delta_interp}
\delta < \frac{\tilde H_{s,\tau}(1)}{m_{\Omega}},
}
where $\tilde H_{s,\tau}$ is the inverse function of  
\bes{
\tilde h_{s,\tau} (z) = \frac{(1+\tau)^s z^{s+1}}{(s+1)!}  \left(1+4z/\pi\right),  \quad z\in(0,\infty),
}
then 
\bes{
A\|f\|^2 \leq \sum_{n\in I}\left \| \cP_n(f) \right \|^2 \leq B \|f\|^2,\quad \forall f\in \rB(\Omega),
}
where $\cP_n(f)$ are  given by \textnormal{\R{interp_operator}} and 
\begin{equation}\label{fusion_fr_bounds}
\begin{aligned}
A &\geq  \left( 1- \frac{(1+\tau)^s(\delta m_{\Omega})^{s+1}}{(s+1)!}  \left(1+\frac{4\delta m_\Omega}{\pi}\right) \right)^2 , \\
B &\leq  \left( 1 + \frac{(1+\tau)^s(\delta m_{\Omega})^{s+1}}{(s+1)!}  \left(1+\frac{4\delta m_\Omega}{\pi}\right) \right)^2.
\end{aligned}
\end{equation}
Equivalently, the set $\left\{ \cP_n  \right\}_{n\in I}$ is a non-orthogonal fusion frame for $\rB(\Omega)$ with weights $v_n = 1$.  
}

\prf{
Let $g(x) = \sum_{n\in I} \cP_n(f)(x)$.  Then
\bes{
\|g\|^2  = \int_{\bbR} \left|\sum_{n\in I} p_n(f)(x) \bbI_{V_{n}}(x)  \right|^2 \D x 
 = \sum_{n\in I} \int_{V_{n}} \left|  p_n(f)(x)  \right|^2 \D x =\sum_{n \in I} \| \cP_n(f) \|^2.
}
Since $f$ is a bandlimited function, it is infinitely continuously differentiable. Also, since for each $n\in I$, $p_n(f)(x)$ is a polynomial of degree at most $s$ that interpolates $f$ at $s + 1$ distinct points $\{x_{n,m}:m=0,\ldots,s\}$ in the closed interval $V_{n}$, a classical result gives that for each $n\in I$ and $x\in V_n$ there exists $\xi_n(x) \in V_{n}$ such that
\be{\label{interp_error}
f(x) - p_n(f)(x) = \frac{f^{(s+1)}(\xi_n(x))}{(s+1)!} \prod_{m=0}^{s} (x - x_{n,m}). 
}
Let $\tilde x_n\in V_{n}$ be such that $|f^{(s+1)}(\tilde x_n)| = \max_{x\in V_{n}}|f^{(s+1)}(x)|$, which again exists because $f$ is bandlimited. Note that, for all $x\in V_{n}$, $\abs{x - x_{n,m}}\leq\delta+h$ for $m\neq0$ and $\abs{x - x_{n,m}}\leq\delta$ for $m=0$. Thus, for all $x\in V_{n}$ we have
\bes{
\left|f(x) - p_n(f)(x)\right| \leq \frac{\left|f^{(s+1)}(\tilde x_n)\right|}{(s+1)!} (1+\tau)^{s}\delta^{s+1}. 
}
Therefore
\eas{
\|f-g\|^2 \! = \! \sum_{n\in I} \int_{V_{n}} \!\! \left|f(x) - p_n(f)(x)\right|^2 \D x \! \leq \! \frac{(1+\tau)^{2s}\delta^{2(s+1)}}{\left((s+1)!\right)^2}  \sum_{n\in I}   \textnormal{meas}(V_{n})  |f^{(s+1)}(\tilde x_n)|^2.
}
By the construction, the points $\{\tilde x_n\}_{n\in I}$ are $2\delta$-dense and $\tilde x_n\in V_{n}$, $n\in I$. Hence, by adapting the proof of Gr\"ochenig's result \cite{GrochenigIrregular} for $s=0$ (to account for the fact that $\tilde{x}_n \neq x_{n,0}$), we get
\bes{
\|f-g\| \leq \frac{(1+\tau)^{s}\delta^{s+1}}{(s+1)!}  \left(1+4\delta m_\Omega/\pi\right) m_{\Omega}^{s+1} \|f\|.
}
The result now follows immediately.
}

\begin{table}
\begin{center}
\scalebox{0.81}{
\begin{tabular}{|c|c|c|c|c|c|c|c|c|c|c|}
\hline
$s$ & 0 & 1 & 2 & 3 & 4 & 5 & 6 & 7 & 8 & 9  \\ \hline 
$\tilde H_{s,1}(1)$  &    0.5766 &  0.7218 &   0.8894 &   1.0626 &  1.2382 &  1.4151 &   1.5928 &  1.7710 &  1.9497 &  2.1287    \\ \hline 
$\tilde H_{s,1/2}(1)$  &  0.5766 &  0.8101 &  1.0458 &  1.2820  & 1.5187 &  1.7558 & 1.9934  & 2.2314 & 2.4696 &  2.7082    \\ \hline 
$\tilde H_{s,1/4}(1)$  &  0.5766  & 0.8710  &  1.1578  &  1.4426  & 1.7270 & 2.0115  & 2.2963  & 2.5815 & 2.8671  &  3.1531   \\ \hline 
$\tilde H_{s,1/8}(1)$  &  0.5766  & 0.9080 & 1.2275  & 1.5440  & 1.8597  & 2.1754  & 2.4914  & 2.8079 & 3.1248  &   3.4422      \\ \hline 
$\tilde H_{s,1/16}(1)$  & 0.5766  &   0.9287 & 1.2669 & 1.6017 & 1.9357  & 2.2696  &  2.6039 & 2.9387 & 3.2740 &   3.6099      \\ \hline
\end{tabular}}
\end{center}\caption{The constant in the bunched sampling density bound \textnormal{\R{delta_interp}}.}
\label{tab:delta_bound_interpol}
\end{table}

The constant $\tilde H_{s,\tau}(1)$ in the density bound obtained by this theorem is explicitly calculated for different values of $s$ and $\tau$ in Table \ref{tab:delta_bound_interpol}. The asymptotic result is given in the following corollary:

\cor{
For large $s$, if 
\bes{
\delta < \frac{1}{(1+\tau)\E}\frac{s+1}{m_{\Omega}},
}
the set $\{ \cP_n \}_{n\in I}$ is a non-orthogonal fusion frame for $\rB(\Omega)$ with bounds as in \textnormal{\R{fusion_fr_bounds}}.
}

\prf{Let $z=\tilde{H}_{s,\tau}(1)$, i.e.\ $\tilde{h}_{s,\tau}(z)=1$. This gives
\bes{
\frac{z}{s+1} \left(1+\tau\right)^{1-\frac{1}{s+1}} \left(1+4z/\pi\right)^{\frac{1}{s+1}} = \frac{((s+1)!)^{\frac{1}{s+1}}}{s+1}.
}
Therefore 
\bes{
\tilde{H}_{s,\tau}(1) \sim \frac{s+1}{(1+\tau)\E}  
}
as $s\rightarrow\infty$.
}

By choosing different form of the interpolation polynomial in \R{interp_operator}, we get different systems. In particular,
for the \textit{Lagrange} form of the interpolation polynomial the operator \R{interp_operator} becomes 
\bes{
\cP_n(f)(x) = \sum_{m=0}^{s}  f(x_{n,m}) L_{n,m}(x) \bbI_{V_n}(x),
}
where $L_{n,m}$ are Lagrange polynomials given by
\be{\label{lag_p_basis}
L_{n,m}(x) = \frac{R_{n,m}(x)}{R_{n,m}(x_{n,m})}, \quad R_{n,m}(x) = \prod_{\substack{0\leq j\leq s\\ j\neq m}} (x-x_{n,j}),
}
and therefore, for the fusion frame operator we have
\bes{
\cS(f)(t) =\sum_{n\in I} \sum_{m=0}^s \sum_{l=0}^{s}  \left( \int_{V_n} L_{n,m}(x) L_{n,l}(x) \D x \right)  f(x_{n,l}) \Phi_{\Omega}(t-x_{n,m}). 
}
On the other hand, if we use the \textit{Newton} form of the interpolation polynomial, we have
\be{\label{operator_newton}
\cP_n(f)(x) = \sum_{m=0}^{s} D_{x_{n,0},\ldots,x_{n,m}}f N_{n,m}(x) \bbI_{V_n}(x),
}
where $D_{x_{n,0},\ldots,x_{n,m}}f$ denotes divided difference of the function $f$ at  $x_{n,0},\ldots,x_{n,m}$ and $N_{n,m}$ is Newton polynomial given by
\be{\label{new_p_basis}
N_{n,m}(x) = \prod_{l=0}^{m-1} (x-x_{n,l}).
}
The fusion frame operator in this case is
\bes{
\cS(f)(t) =\sum_{n\in I} \sum_{m=0}^s \sum_{l=0}^{s}  \left( \int_{V_n} N_{n,m}(x) N_{n,l}(x) \D x \right)  D_{x_{n,0},\ldots,x_{n,l}} f D_{x_{n,0},\ldots,x_{n,m}}\Phi_{\Omega}(t-\cdot). 
}
Moreover, this approach allows us to consider the following more general setting. Suppose that we are additionally given $k$ derivatives at the points of the bunched sampling set $\{x_{n,m}\}_{n\in I,0\leq m\leq s}$, i.e.\ the given data is
\bes{
f^{(j)}(x_{n,m}), \quad n\in I, \ m=0,\ldots,s, \ j=0,\ldots,k.
}
Now, for each $n\in I$, we can define the unique interpolation polynomial $p_n(f)$ such that
\bes{
p_n^{(j)}(f)(x_{n,m}) = f^{(j)}(x_{n,m}), \quad m=0,\ldots,s,  \ j=0,\ldots,k.
}
In this case, we can use the \textit{Hermite} form of the interpolation polynomial and set
\bes{
\cP_n(f)(x) = \sum_{j=0}^{k} \sum_{m=0}^{s} f^{(j)}(x_{n,m}) c_{n,m,j}(x) \bbI_{V_n}(x),
}
where
\bes{
c_{n,m,j}(x) = L_{n,m}^{k+1}(x) \frac{(x-x_{n,m})^j}{j!} \sum_{i=0}^{k-j} \frac{(x-x_{n,m})^i}{i!} R_{n,m}^{k+1}(x_{n,m}) \frac{d^i}{dx^i} R_{n,m}^{-(k+1)}(x_{n,m}),
}
and $L_{n,m}$, $R_{n,m}$ are as in \R{lag_p_basis}, see \cite{TraubHermite}. Since \R{interp_error} now reads as
\bes{
f(x) - \cP_n(f)(x) = \frac{f^{((s+1)(k+1))}(\xi_n(x))}{((s+1)(k+1))!} \prod_{m=0}^{s} (x - x_{n,m})^{k+1}, 
}
we obtain an additional $k+1$ factor in the density bound, i.e.\ the density condition now reads
\bes{
\frac{(1+\tau)^{s(k+1)}(\delta m_{\Omega})^{(s+1)(k+1)}}{(s+1)!(k+1)!}  \left(1+4\delta m_\Omega/\pi\right)  < 1,
}
which for large $s$ and large $k$ leads to
\bes{
\delta < \frac{1}{(1+\tau)\E} \frac{(s+1)(k+1)}{m_{\Omega}}.
}
Thus a combination of bunched and derivative sampling increases the maximal allowed density by a multiplicative factor in $s+1$ (number of bunched points) and $k+1$ (number of derivatives).

\subsection{Bunched sampling and frames}\label{ss:bunched_frames}

It transpires that the use of the Newton form of the interpolating polynomial also allows one to relate bunched sampling to a frame, as opposed to a fusion frame.  Let us define $\cP_n$ as in \R{operator_newton}.
Since $D_{x_{n,0},\ldots,x_{n,m}} f$ is just a linear combination of the function $f$ evaluated at the points $x_{n,0},\ldots,x_{n,m}$ and since $f(x)=\ip{f(t)}{\Phi_{\Omega}(t-x)}$ with $\Phi_{\Omega}$ defined by \R{Phi_fun}, we can write 
\be{\label{little_phi}
D_{x_{n,0},\ldots,x_{n,m}} f = \ip{f}{\phi_{n,m}},\quad \phi_{n,m}(t) = D_{x_{n,0},\ldots,x_{n,m}}\Phi_{\Omega}(t-\cdot).
}
We now have the following:

\prop{\label{p:bunched_samp_finitediff}
Suppose that $\{x_{n,m}\}_{n\in I,0\leq m\leq s}\subseteq\bbR$ is the bunched sampling set with density $\delta$ and bunch width $h=\tau\delta$, where $\tau\in(0,1]$. Let $\{V_n\}_{n\in I}$ be  the Voronoi regions corresponding to the points $\{x_{n,0}\}_{n\in I}$. If
\bes{
\delta < \frac{\tilde H_{s,\tau}(1)}{m_{\Omega}},
}
where $\tilde H_{s,\tau}$ is as in Theorem \textnormal{\ref{t:bunched_samp}},  
then 
\be{\label{finitediff_frame}
A\|f\|^2 \leq \sum_{n\in I} \sum_{m=0}^{s}  \mu_{n,m}  \left|  D_{x_{n,0},\ldots,x_{n,m}} f \right|^2  \leq B \|f\|^2, \quad \forall f\in \rB(\Omega),
}
where $\mu_{n,m}=m!\int_{V_n} \left|N_{n,m}(x)\right|^2 \D x$, $N_{n,m}$ are  given by \textnormal{\R{new_p_basis}} and 
\be{\label{fr_bounds_divided_dif}
A \geq \! \frac{1}{\E} \! \left(\! 1 \!- \! \frac{(1+\tau)^s(\delta m_{\Omega})^{s+1}}{(s+1)!} \! \left(\!1\!+\!\frac{4\delta m_\Omega}{\pi}\right)\!\! \right)^2 \!, \ B\! \leq\!  \frac{\left(1+ 2(1+\tau)\delta m_{\Omega}/\pi \right)^2 \E^{((1+\tau)\delta m_{\Omega})^2}}{(1+\tau)^2}\!.
} 
Equivalently, if $\phi_{n,m}$ is as in \textnormal{\R{little_phi}}, the set $\left\{ \sqrt{\mu_{n,m}} \phi_{n,m} :  n\in I, m=0,\ldots,s \right\}$
is a frame for $\rB(\Omega)$. 
}
\prf{
As before, let
\bes{
g(x) = \sum_{n\in I} \sum_{m=0}^{s} D_{x_{n,0},\ldots,x_{n,m}} f  N_{n,m}(x) \bbI_{V_n}(x).
}
Now we have
\eas{
\|g\|^2  &=  \sum_{n\in I} \int_{V_{n}} \left|\sum_{m=0}^{s}  D_{x_{n,0},\ldots,x_{n,m}} f  N_{n,m}(x)  \right|^2 \D x\\
&\leq \sum_{m=0}^{s} \frac{1}{m!} \sum_{n\in I} \sum_{m=0}^{s} m! \left( \int_{V_n} \left|N_{n,m}(x)\right|^2 \D x \right) \left| D_{x_{n,0},\ldots,x_{n,m}} f  \right|^2
}
and hence
\bes{
\sum_{n\in I}  \sum_{m=0}^{s}  \mu_{n,m}  \left| D_{x_{n,0},\ldots,x_{n,m}} f   \right|^2 \geq  \E^{-1} \left( \|f\| -  \|f - g\| \right)^2.
}
In the proof of Theorem \ref{t:bunched_samp} we obtained 
\bes{
\|f-g\| \leq \frac{(1+\tau)^s(\delta m_{\Omega})^{s+1}}{(s+1)!}  \left(1+4\delta m_\Omega/\pi\right) \|f\|,
}
and therefore for the lower frame bound we get
\bes{
A\geq \E^{-1} \left( 1- \frac{(1+\tau)^s(\delta m_{\Omega})^{s+1}}{(s+1)!}  \left(1+4\delta m_\Omega/\pi\right) \right)^2.
}
For the upper frame bound first note that
\bes{
\mu_{n,m} = m! \int_{V_n} \left|\prod_{l=0}^{m-1} (x-x_{n,l})\right|^2 \D x \leq m! (1+\tau)^{2(m-1)}\delta^{2m} \text{meas}(V_n).
}
Also, by the mean value theorem for  divided differences, for every $n\in I$ there exists $\tilde x_{n}\in [x_{n,0}-h,x_{n,0}+h]$ such that
\bes{
D_{x_{n,0},\ldots,x_{n,m}} f  = \frac{f^{(m)}(\tilde x_{n})}{m!}.
}
Now, as in the proof of Theorem \ref{t:main_fram_1D}, since the points $\{\tilde x_n\}_{n\in I}$ are $(\delta+\tau\delta)$-dense, we obtain
\eas{
\sum_{n\in I} \sum_{m=0}^{s}  \mu_{n,m}  \left| D_{x_{n,0},\ldots,x_{n,m}} f   \right|^2 &\leq \frac{1}{(1+\tau)^2} \sum_{m=0}^{s} \frac{((1+\tau)\delta)^{2m}}{m!} \sum_{n\in I}  \text{meas}(V_n)   \left|  f^{(m)}(\tilde x_{n})  \right|^2\\
&\leq \frac{\left(1+ 2(1+\tau)\delta m_{\Omega}/\pi \right)^2 \E^{((1+\tau)\delta m_{\Omega})^2}}{(1+\tau)^2}  \|f\|^2,
}
and the estimate for the upper frame bound follows.
}

In the limit, when the bunch width $h$ becomes very small and the number of bunched points $s$ very large, from this proposition we obtain precisely the one-dimensional derivative result given in Theorem \ref{t:main_fram_1D} for large number of derivatives $k$:

\cor{For large $s$ and small $\tau$, if 
\bes{
\delta < \frac{1}{\E}\frac{s+1}{m_{\Omega}},
}
then $\left\{ \sqrt{\mu_{n,m}}  \frac{d^{m}}{dx^m}\Phi_{\Omega}(\cdot-x_{n,0})  : \mu_{n,m}=\frac{1}{m!}\int_{V_n}(x-x_{n,0})^{2m}\D x, n\in I, m=0,\ldots,s \right\}$ is a frame for $\rB(\Omega)$ with the frame bounds satisfying \R{fr_bounds_divided_dif}.
}

\prf{
Consider the sum \R{finitediff_frame} as $\tau\rightarrow 0$. For $x_{n,0},\ldots,x_{n,m} \!\in\! [x_{n,0}-\tau\delta,x_{n,0}+\tau\delta]$  
\eas{
&\lim_{\tau\rightarrow0} \sum_{n\in I} \sum_{m=0}^{s}  m! \int_{V_n}\abs{N_{n,m}(x)}^2\D x  \left|  D_{x_{n,0},\ldots,x_{n,m}} f \right|^2 \\
&=  \sum_{n\in I} \sum_{m=0}^{s} \frac{1}{m!}\int_{V_n}(x-x_{n,0})^{2m}\D x |f^{(m)}(x_{n,0})|^2.
}
This holds due to dominated convergence theorem, since for any $\tau$, $n$ and $m$ 
\bes{
m! \int_{V_n}\abs{N_{n,m}(x)}^2\D x  \left|  D_{x_{n,0},\ldots,x_{n,m}} f \right|^2 \leq  \text{meas}(V_n) (2\delta)^{2m} |f^{(m)}(\tilde x_n)|^2,
}
where $\tilde x_n\in V_{n}$ is such that $|f^{(m)}(\tilde x_n)| = \max_{x\in V_{n}}|f^{(m)}(x)|$.

For the density condition, as before, let $z=\tilde{H}_{s,\tau}(1)$. Since $1+\tau\sim1$ as $\tau\rightarrow 0$, this gives
\bes{
\frac{z}{s+1}  \left(1+4z/\pi\right)^{\frac{1}{s+1}} \sim \frac{((s+1)!)^{\frac{1}{s+1}}}{s+1},\quad \tau\rightarrow 0,
}
and hence $\tilde{H}_{s,\tau}(1) \sim (s+1)/\E$  as $\tau\rightarrow 0$ and $s\rightarrow\infty$. 
}

Therefore, for the large number of bunched sampling points $s$ such that the width of all bunches is small, we obtain the same result as when sampling $s$ derivatives.  

\section{Conclusions}
In this paper, we present several density bounds as sufficient guarantees for stable recovery of a bandlimited function from measurements of it and its first $k$ derivatives. We also have proved the linear growth of $\delta$-density with $k+1$. However, the constant of proportionality $1/\E$ is rather small compared to the case of equispaced samples where the corresponding constant is $\pi/2$. Therefore, it would be of interest to see how these bounds can be improved in both the univariate and multivariate case. We also believe that the perturbation theorem given in this paper can be improved. Moreover, for the sake of better understanding of the bound on the perturbation distance $\epsilon$, it is important to analyse the behaviour of the ratio $A/B$ with $k$. This is left for future work. 

As we have seen, a related problem to derivatives sampling is so-called bunched sampling. This sampling strategy also leads to increased $\delta$-bound and, asymptotically, it approximates the derivatives sampling. Much as in the derivatives case, it remains open to improve this density bound. Also, it would be important to generalize these results to the multivariate case and therefore broaden the range of their applications. Let us note that in higher dimensions, well-posedness of the bunched points and the possibility of constructing an unique multivariate interpolation polynomial complicates dramatically. Therefore, we leave this problem for future investigations.

One might notice that in this paper we analyse two examples---derivatives and bunched sampling---both appearing at the end of Papoulis' paper \cite{PapoulisGenSamp}. Although these examples are of interest in applications by themselves, the remaining  problem is to analyse a general setting given in Papoulis' paper in the context of nonuniform sampling. It remains open to see what happens with the sampling density when instead of $H_{\alpha}(\omega)=(-\I\omega)^{\alpha}$ one has more general functions $H_{\alpha}$ and a nonuniform set of sampling points.

\section*{Acknowledgements}
The authors would like to thank Akram Aldroubi, Karlheinz Gr\"ochenig, Maarten Van De Hoop and Ilya Krishtal for useful discussions. Additionally, the authors would like to thank to the participants of the ICERM Research Cluster ``Computational Challenges in Sparse and Redundant Representations'' for providing a stimulating and interactive research environment.

BA was supported by the NSF DMS grant 1318894. MG and AH were supported by the EPSRC grant  EP/N014588/1 for the EPSRC Centre for Mathematical and Statistical Analysis of Multimodal Clinical Imaging. AH was also supported by a Royal Society University Research Fellowship as well as the EPSRC grant EP/L003457/1.

\bibliographystyle{siam}
\bibliography{DerivSampRefs}

\end{document}